\documentclass[10pt]{article}
\usepackage{latexsym,amssymb,amsmath,url,authblk}

\def\N{\mathbb{N}}
\def\Q{\mathbb{Q}}

\def\R{\mathbb{R}}

\def\proof{\par\medskip\noindent{\em Proof. }}
\def\eproof{\hfill{$\Box$}\bigskip}

\def\tec{\hspace{-1.6mm}{\bf. }}
\def\ds{\dots}
\def\sus{\subset}
\def\al{\alpha}
\def\be{\beta}
\def\ga{\gamma}

\def\ep{\varepsilon}
\def\cc{\colon}

\newtheorem{thm}{Theorem}[section]
\newtheorem{prop}[thm]{Proposition}

\newtheorem{lem}[thm]{Lemma}

\newtheorem{defi}[thm]{Definition}
\newtheorem{exam}[thm]{Example}

\title{Around Hilbert's theorem: the center of a~circle is not constructible by straightedge alone}
\author{Martin Klazar}
\affil{Department of Applied Mathematics, Faculty of Mathematics and Physics, Charles University, Malostransk\'e n\'am\v est\'\i\ 25, 118 00 Praha 1, Czechia\\ e-mail: {\tt klazar@kam.mff.cuni.cz}}
\date{}

\begin{document}
\maketitle

\begin{abstract}
In order to state the theorem in the title formally and to review its rigorous proof, we extend and make more precise the Uspenskiy--Shen--Akopyan--Fedorov model of Euclidean constructions with arbitrary points; we also introduce formalizations for  infinite configurations and for the projective plane. We exemplify the proof method by
simpler and not so well known results that it is impossible to construct the unit length, or a~given point, by compass and straightedge from nothing by means of classical arbitrary points. On the other hand we construct any given point by compass and straightedge from nothing by means of arbitrary points 
determined by horizontal segments. We quote a~``proof'' of Hilbert's theorem from the literature and explain why it is problematic. 
We rigorously prove Hilbert's theorem and present three variants of it, the last one for the projective plane.

\end{abstract}

\section{Introduction}\label{intro}

In 1913, D.~Cauer\footnote{By footnote 137 in D.\,E.~Rowe \cite[p. 229]{rowe}, Detlef Cauer (1889--1918) 
was a~son of the classical philologist Paul Cauer, studied mathematics in Kiel, Berlin, M\"unster and G\"ottingen, later was an 
assistant of E. Landau, and in April 1918 fell in WWI in Belgium. For some more information on him see \cite{wiki_caue} .} 
mentioned  in \cite{caue} that D.~Hilbert had proved during his lectures that it is impossible to construct the center of a~given circle only by straightedge; 
D.~Hilbert did not publish his proof. 
D.~Cauer generalized Hilbert's argument and proved that it is impossible to construct only by straightedge the centers of two given 
circles, if the circles are disjoint and not concentric (i.e., do not have common center). For two intersecting or two concentric circles 
straightedge-only constructions of the centers are known, see \cite[p. 96]{akop_fedo} or \cite[p. 93]{caue} for both, and \cite[p. 173]{rade_toep} 
for a~detailed explanation of the former construction. Forty years later, C.~Gram \cite{gram} found an error in Cauer's proof: it works equally well when 
besides the two circles also a point on the line connecting their centers is given, but then C.~Gram gave a straightedge-only construction of the 
centers. In fact, Cauer's theorem is wrong. A.~Akopyan and R.~Fedorov \cite[p. 97/8]{akop_fedo} gave straightedge-only constructions of the centers 
for certain pairs of disjoint and non-concentric circles. For example, if a~circle $k_1$ lies inside a~circle $k_2$ and there exists a~quadrilateral 
inscribed in $k_2$ and circumscribed around $k_1$, then one can construct the centers of $k_1$ and $k_2$ only by straightedge. A.~Akopyan 
and R.~Fedorov could save from Cauer's theorem the result \cite[Theorem 1.1]{akop_fedo}: ``There exist two circles whose centers cannot be constructed using
only a straightedge.''

Gy.~Strommer \cite{stro} published another strengthening of Hilbert's theorem: even if besides the circle two perpendicular lines (not crossing 
exactly in the center of the circle), or a~line (not going through the center) with three points on it marking two segments with equal lengths, 
are given, one still cannot construct the center only by straightedge (\cite[part 3 from p. 97]{stro}). Unfortunately, his proof rests on the 
same transformation fallacy, discussed below, as other published proofs of Hilbert's theorem. Gy.~Strommer described simple and explicit 
deforming transformations which we review in Proposition~\ref{stro_map} and use first in the correct proof of Hilbert's theorem in Theorem~\ref{hilb_thm}, and then in 
Theorem~\ref{concr_hilb_thm} in a~concrete deterministic Hilbert's theorem. 

During the 20th century Hilbert's theorem and its proof were mentioned in several expository books: in R.~Courant and H.~Robbins \cite{cour_robb}, in M.~Kac 
and S.~M.~Ulam \cite{kac_ulam}, in H.~Rademacher and O.~Toeplitz \cite{rade_toep}, and in some others; see A.~Shen \cite{shen} and V.~Uspenskiy 
and A.~Shen \cite{shen_uspe} for more references. Any attempted construction of the center of a circle only by straightedge, or even when 
compass is allowed, that starts with the bare circle is non-deterministic. To begin one has to select an arbitrary point in the 
plane. In fact, to get anywhere one has to select at least three distinct arbitrary points. But isn't one arbitrary point enough? Why not to pick 
exactly the center of the given circle? This gives a~very short construction with just one step. Yes, we are cheating, but why? Until recently the problem with Hilbert's 
and Cauer's proofs and their variants was that they did not use any precise definition of Euclidean constructions with arbitrary 
points, a~definition that explains why selecting the sought-for center as an arbitrary point and ``refuting'' by this Hilbert's theorem is 
not allowed. See \cite{shen_uspe} for the history of attempts to deal with arbitrary points in Euclidean constructions. When finally in 2017/18, with 
a~century delay, precise definitions of Euclidean constructions with arbitrary points were proposed by A.~Akopyan and R.~Fedorov \cite{akop_fedo}, 
A.~Shen \cite{shen}, and V.~Uspenskiy and A.~Shen \cite{shen_uspe}, it became clear that Hilbert's argument is fallacious and does not prove the result. 

We want to report to the reader on this interesting, even if somewhat embarrassing, development. In Section~2 we quote from the literature one
``proof'' of Hilbert's theorem and explain why it is not sufficient. It is not so strange that D.~Hilbert and others erred, every mathematician 
knows the terrible power of wishful thinking\,---\,how easily one gets convinced  that the current plausible argument already 
is the desired rigorous proof, or that it {\em could} be easily made in one by filling in just few inessential technical details. It is more disconcerting 
and worrying that Hilbert's ``proof''  was uncritically taken over both in popular accounts and research articles, and that it took over 100 years to  recognize it clearly as fallacious. The first correct proof of Hilbert's theorem was given by A.~Akopyan and R.~Fedorov in \cite{akop_fedo} in 2017.

In Section~2 we present in Theorem~\ref{hilb_thm} a~rigorous proof of Hilbert's theorem, with more details than in 
\cite{akop_fedo} and \cite{shen}. But we begin the section with Definition~\ref{defi_EC} of $\mathrm{EC}(\mathcal{S})$, Euclidean construction with $\mathcal{S}$-arbitrary points, one of the main results of our article. It is similar to the 
{\em game definition} in A.~Shen \cite{shen} and in V.~Uspenskiy and A.~Shen \cite{shen_uspe}, but it is more general since we allow any set system 
${\cal S}$ for determining arbitrary points, not just open sets, and it is also more formal and precise since we coach it in terms of rooted trees as a~concrete and precise
set-theoretic structure. In our approach an Euclidean construction with arbitrary points is a~concrete and ``tangible'' set-theoretic object; such 
concreteness and rigour is still missing in the approaches of \cite{akop_fedo,shen,shen_uspe}. Then we give Example~\ref{example} of our definitions and illustrate by Example~\ref{ex_schreiber} the approach to Euclidean constructions by P.~Schreiber \cite{schr}. In 
Propositions~\ref{prop_constr} and \ref{equi_mode} we establish general properties of our model and after auxiliary Lemma~\ref{e_closed} we illustrate it 
by proving in Theorem~\ref{impo_unit} a~simple and unjustly unknown impossibility result: there is no Euclidean construction with {\em classical} 
arbitrary points (determined by open sets) that
uses compass and straightedge, starts from the empty configuration, and constructs two points with distance $1$. In Proposition~\ref{po_unit} we show 
that such construction is possible with $\mathcal{U}$-arbitrary points that are determined by horizontal segments. Similarly, 
in Proposition~\ref{hilb_inval} we present an Euclidean construction of the center of a~given circle only by straightedge, but with the help of 
$\mathcal{U}$-arbitrary points. This does not refute Hilbert's theorem because stronger arbitrary points than the classical ones are  used, but it convincingly refutes any ``proof'' of the theorem that lacks 
precise specification of arbitrary points. We give an example of such a~proof in a~verbatim quote, and its translation, of a~passage from \cite{rade_toep} and discuss its shortcomings.  Theorem~\ref{hilb_thm} rigorously states and proves
Hilbert's theorem. For the proof we need certain deformation maps which Gy.~Strommer described conveniently for us in 
\cite{stro}. We adapt his construction in Propositions~\ref{stro_map} and \ref{map_fp}. 
We deliberately eliminate from the proof of Theorem~\ref{hilb_thm} any projective element (they are used in other proofs of Hilbert's theorem) so that they cannot hide any error. Theorems~\ref{imp_point} and \ref{po_point} are devoted 
to the simplest-to-state problem in Euclidean constructions with arbitrary points: construct a~given point, 
say the origin, by compass and straightedge from nothing. This may sound as a~trivial problem but in reality is not. In the former theorem we show\,---\,we omit the proof as it is very similar 
to that of Theorem~\ref{impo_unit}\,---\,that the construction does not exist if only classical arbitrary points are allowed. The latter theorem however shows that the 
construction is possible by means of $\mathcal{U}$-arbitrary points.
 
Section 3 contains three variants of Hilbert's theorem. 1.~In Theorem~\ref{hilb_thm_det_strong},  stated in Theorem~\ref{hilb_thm_det_strong3} for infinite configurations, we show that any circle $k$ in the plane 
has a~countable and dense subset $Y\sus k$ such that for no finite subset $Z\sus Y$ there is an Euclidean construction starting from $k$ and the point set $Z$ 
that deterministicly constructs the center of $k$  only by straightedge. 2.~In Theorem~\ref{concr_hilb_thm} we present, using transcendence 
of the numbers $\sin1$ and $\cos1$, an explicit example of such set 
$Y$ for a~particular circle. 3.~In Theorem~\ref{proj_hilb_thm} we give a precise statement and proof of Hilbert's theorem in the  projective plane; 
the proof rests on the projective version of Strommer's map.
To our knowledge this is the first rigorous treatment of a~projective version of Hilbert's theorem. We hope to continue our investigation of Euclidean constructions with $\mathcal{S}$-arbitrary points 
in \cite{klaz}.

\section{Euclidean constructions with $\mathcal{S}$-arbitrary\\ points and Hilbert's theorem}

Let $\N=\{1,2,\ds\}$ be the natural numbers, $\N_0=\N\cup\{0\}=\{0,1,\ds\}$ be the nonnegative integers, and 
$\omega=\{0,1,\ds\}=\{\emptyset,\{\emptyset\},\ds\}$ be the first infinite ordinal. By $\R$ (resp. $\Q$) we denote the real (resp. rational) numbers. 
For a~set $A$, called an {\em alphabet}, 
we consider sequences $u=(a_i)$ with entries $a_i\in A$. We call the $a_i$ the {\em letters
of $u$}. If $u$ is finite, we call it a~{\em word (over $A$)} and write it as $u=a_1a_2\ds a_m$ with $m\in\N_0$. For $m=0$ we get the 
{\em empty word} $u=\emptyset$. We denote the set of words over $A$ by $A^*$. We will consider also certain infinite sequences $u$ 
with entries in $A$, called {\em infinite words (over $A$)}.
We write them as $u=a_0a_1\ds$ if $u$ is indexed by the elements of $\omega$ (for example, $u$ is an infinite path in a~rooted tree), or as 
$u=a_0a_1\ds;a_{\omega}a_{\omega+1}\ds a_{\omega+i}$, where $i\in\N_0$, if $u$ is indexed by the elements of the ordinal $\omega+i+1$
(see the comment after Theorem~\ref{hilb_thm_det_strong}). For a~possibly infinite set $X$ we denote by $|X|$ its cardinality. 

A~\emph{rooted tree $T=(r,V,E)$} is a~triple of a {\em root $r\in V$}, 
a~set $V$ of {\em vertices}, and a set $E\sus V\times V$ of {\em edges} that satisfies the following condition. For every vertex $u\in V$ there 
is a~unique {\em walk from $r$ to $u$}: a~unique word $u_1u_2\ds u_k\in V^*$, $k\in\N$, such that $u_1=r$, $u_k=u$, and
$(u_i,u_{i+1})\in E$ for every $i=1,2,\ds,k-1$. The uniqueness implies that each of these walks is in fact a~{\em path}, $u_i\ne u_j$ for $i\ne j$. 
More generally, a~{\em walk in $T$} is any word $u_1u_2\ds u_m\in V^*$, or any infinite word $u_0u_1\ds$ with $u_i\in V$, such that $(u_i,u_{i+1})\in E$ 
for every $i$. By the uniqueness each walk is a~path, no vertex is repeated. From the uniqueness it also follows  that every vertex $u\ne r$ has 
in $T$ in-degree $1$, and that $r$ has in-degree $0$, where the {\em in-degree} 
of a~vertex $u$ is the number of vertices $v$ with $(v,u)\in E$. Similarly, the {\em out-degree of $u\in V$} is the number of vertices $v$ with 
$(u,v)\in E$; the vertices $v$ are the {\em children of $u$} and $u$ is their {\em parent}. Out-degrees may attain any value, and the vertices with 
out-degree $0$ are called {\em leaves of $T$}. If a~walk in $T$ is maximal, cannot be prolonged in either way, it is called a~{\em branch in $T$}. 
Every branch in $T$ starts in $r$ and either finishes in a~leaf and is finite, or continues forever and is infinite. 

We denote by ${\cal P}=\R^2$ the set of {\em points}, by ${\cal L}$ the set of {\em lines} in $\R^2$, and by ${\cal C}$ the set of {\em circles} with positive radii 
in the affine plane $\R^2$. We call the elements of ${\cal L}\cup{\cal C}$ {\em curves}. For $a,b\in{\cal P}$ with $a\ne b$ we denote by 
$l(a,b)\in{\cal L}$ the line going through the points $a$ and $b$, by $ab$ the segment spanned by them, and by $|ab|$ its length. For $a,b,c\in{\cal P}$ with 
$b\ne c$ we denote by $k(a,b,c)\in{\cal C}$ the circle with 
center $a$ and radius $|bc|$. For $a,b\in{\cal P}$ we set $k(a,b,b)=a\in{\cal P}$; these degenerated circles are important 
as they enable us to repeat any selected point. In the case of two non-parallel distinct lines $\kappa$ and $\ell$ we abuse set notation 
and write $\kappa\cap\ell=p\in{\cal P}$ for their intersection point $p$, instead of the correct $\kappa\cap\ell=\{p\}$. If the lines $\kappa$ and 
$\ell$ are parallel, i.e. $\kappa\cap\ell=\emptyset$, we write $\kappa\parallel\ell$.

Let ${\cal S}$ be a~possibly empty set of nonempty subsets of the affine plane $\R^2$ (in Section~3 we work also with the projective plane $\mathbb{P}_2$); it is the set of possible locations of arbitrary points. 
For technical reasons we assume that $({\cal P}\cup{\cal L}\cup{\cal C})\cap{\cal S}=\emptyset$. If we want to work, for example, with ${\cal S}={\cal C}$, 
we set ${\cal S}={\cal C}\times\{0\}$ and modify accordingly all definitions. The set system ${\cal S}$ may be the system ${\cal O}$ 
of all nonempty {\em open subsets of $\R^2$} (in the Euclidean topology), or the system ${\cal D}$ of all {\em open discs $D$ in $\R^2$ with positive radii}, or the {\em empty 
system ${\cal S}=\emptyset$} leading to deterministic constructions, or the system ${\cal J}=\{\{s,t\}\;|\;s,t\in\R^2, s\ne t \}$ of all {\em two-element point sets}, 
or the system 
$$
{\cal U}=\{[a,\,b]\times\{c\}\;|\;a,\,b,\,c\in\R,\,a<b\}
$$ 
of all {\em proper horizontal segments}, or some other set system. 

\begin{defi}[$\mathrm{EC}(\mathcal{S})$]\tec\label{defi_EC}
An {\em Euclidean construction with ${\cal S}$-arbitrary points}, abbreviated $\mathrm{EC}(\mathcal{S})$, is any rooted tree $T=(r,V,E)$ with the next described structure. Its vertices 
$$
u=a_1a_2\ds a_m\in V\sus({\cal P}\cup{\cal L}\cup{\cal C}\cup{\cal S})^*,\ m\in\N_0\;,
$$ 
are words over the alphabet of points, lines, circles, and 
elements of ${\cal S}$. Each vertex $u\in V$ is either {\em deterministic} with $m\ge0$, $a_m\not\in{\cal S}$ and out-degree $1$ or $0$, or 
{\em non-deterministic} with $m\ge1$, $a_m\in{\cal S}$ and $|a_m|$ children. Always $r\in({\cal P}\cup{\cal L}\cup{\cal C})^*$. Note that every leaf in $T$ 
is a~deterministic vertex. The children of any vertex $u$ (that is displayed above) are determined by exactly one of the following six rules. In the first 
five the vertex $u$ is assumed deterministic.

\begin{enumerate}
    \item Construction ends: $u$ is a leaf with no child.
    \item A~new line: $u$ has one child $v=a_1a_2\ds a_{m+1}$ with $a_{m+1}=l(a_i,a_j)\in{\cal L}$ for some indices $1\le i<j\le m$ such that 
          $a_i,a_j\in{\cal P}$ and $a_i\ne a_j$.
    \item A~new circle or repeated point: $u$ has one child $v=a_1a_2\ds a_{m+1}$ with $a_{m+1}=k(a_i,a_j,a_k)\in{\cal C}\cup{\cal P}$ for some indices $1\le i,j,k\le m$ 
          such that $a_i,a_j,a_k\in{\cal P}$.
    \item A~new intersection point: $u$ has one child $v=a_1a_2\ds a_{m+1}$ with $a_{m+1}\in{\cal P}$ being an intersection point of two curves $a_i\ne a_j$, 
          $1\le i<j\le m$, in $u$. 
    \item A~new location for arbitrary points: $u$ has one child $v=a_1a_2\ds a_{m+1}$ with $a_{m+1}\in{\cal S}$. 
    \item New arbitrary points: $u$ is non-deterministic with $a_m\in{\cal S}$ and has $|a_m|$ children $v=a_1a_2\ds a_{m+1}$, one for each point $a_{m+1}\in a_m$.
\end{enumerate}
\end{defi}
This is the definition of an $\mathrm{EC}(\mathcal{S})$, more precisely of its finitary affine form, which is the form we mostly use here. Later we will briefly consider also the version $\mathrm{EC}_{\infty}(\mathcal{S})$ for infinite countable configurations, and at greater length the version
$\mathrm{EC}_{\mathrm{pr}}(\mathcal{S}_{\mathrm{pr}})$ for the projective plane. The main innovation is that locations $S\in\mathcal{S}$ for arbitrary points are treated on par with points, lines and circles.

We say that an $\mathrm{EC}(\mathcal{S})$ $T$ is a~{\em straightedge construction} if it has no edge $(u,v)$ obtained in the case $a_j\ne a_k$ of rule 
$3$, i.e. compass is forbidden. Similarly, $T$ is a~{\em compass  construction} if it has no edge $(u,v)$ obtained by rule $2$, i.e. straightedge is forbidden. 
In a~{\em general construction $T$} both devices are allowed. We call the three previous possibilities the {\em types of $T$}. We say that $T$ is {\em deterministic} if 
${\cal S}=\emptyset$; equivalently, $T$ has no edge $(u,v)$ obtained by rules $5$ and $6$. Then  for every $u\in V$ one has that
$u\in({\cal P}\cup{\cal L}\cup{\cal C})^*$ and $T$ consists of a single branch. We say that $T$ is {\em terminating} if all branches in $T$ are finite.  
$T$ is {\em classical} if ${\cal S}={\cal D}$, the system of open discs. If
$$
{\cal K}\sus({\cal P}\cup{\cal L}\cup{\cal C})^*\;\text{ and }\;r\in({\cal P}\cup{\cal L}\cup{\cal C})^*
$$
satisfy $u\cap r=\emptyset$ for every $u\in{\cal K}$ (i.e., no letter in $u$ appears as a~letter in $r$), we say that ${\cal K}$ and $r$ are {\em separated}. 
The words in ${\cal K}$ are the allowed terminal configurations and the word $r$ is the initial configuration. We say that an {\em $\mathrm{EC}(\mathcal{S})$ of a~certain 
type constructs ${\cal K}$  from $r$} if there exists a~terminating $\mathrm{EC}(\mathcal{S})$ $T=(r,V,E)$ of the stated type, with the prescribed root $r$,
and such that every leaf $u=a_1a_2\ds a_m$ in $T$ has a~final segment $a_ja_{j+1}\ds a_m\in{\cal K}$. We say that an {\em $\mathrm{EC}(\mathcal{S})$ of a~certain 
type weakly constructs ${\cal K}$ from $r$} if there exists a~terminating $\mathrm{EC}(\mathcal{S})$ $T=(r,V,E)$ of the stated type, with the prescribed root $r$, and such that every 
leaf $u=a_1a_2\ds a_m$ in $T$ has a~(not necessarily contiguous) subsequence that is a~permutation of a~word in ${\cal K}$. 

We remark at this point that the principle of a~successful non-deterministic geometric construction that every possible way of performing it results in the 
desired object, appears in a~form already in Yu.~Manin \cite{mani}, as quoted at p.~13  of the arXiv version of \cite{shen_uspe}.

We allow {\em non-uniform constructions}, which means that the kinds of steps after selecting an arbitrary point $p\in S\in{\cal S}$ may 
depend on $p$. For example, for an arbitrary point $p\in S$ some later step may be drawing a~line through two already constructed points, for another 
arbitrary point $q\in S$, $q\ne p$, that step may be selecting an $S'\in{\cal S}$, and so on. We look at uniformity of Euclidean constructions in more detail in \cite{klaz}. The next example illustrates the above 
notions and is a~uniform construction in which kinds of steps do not depend on selected arbitrary points. The same holds for the constructions in 
Propositions~\ref{po_unit} and \ref{hilb_inval}. In contrast, the construction in Theorem~\ref{po_point} is non-uniform.

\begin{exam}[equilateral triangles]\tec\label{example}
Let
$$
{\cal K}=\{abc\in{\cal P}^3\;|\;|ab|=|ac|=|bc|>0\}\ \mathrm{and}\ r=\emptyset\;.
$$
There exists a~compass $\mathrm{EC}(\mathcal{D})$ $T$ that constructs ${\cal K}$, an equilateral triangle, from nothing. It is 
clear that $T$ is classical and terminating.
\end{exam}
We describe this construction $T=(\emptyset,V,E)$. It consists of $\mathfrak{c}$ (continuum many) branches, each of which has $10$ vertices. The second and fourth vertex (counted from the root)
have out-degree $\mathfrak{c}$ but other vertices are deterministic. Each of the $\mathfrak{c}$ 
leaves $u$ of $T$ has the same form
$$
u=D_1\,p_1\,D_2\,p_2\,C_1\,C_2\,p_3\,p_1\,p_2\;,
$$
where $D_1$ and $D_2$ are common to all $u$ but the other seven letters depend on $u$. The $D_i\in{\cal D}$, $i=1,2$, are open discs with radii 
$1$ and respective centers $(0,0)$ and $(0,3)$ (rule $5$), $p_i\in D_i$ are two arbitrary points (rule $6$), $C_i=k(p_i,p_1,p_2)\in{\cal C}$ are two circles 
with centers in the two points and radii equal to their distance (rule $3$), $p_3\in C_1\cap C_2$ is one of the two intersection points of the two circles 
(rule $4$), and the final points $p_1=k(p_1,p_1,p_1)$ and $p_2=k(p_2,p_2,p_2)$ are repeated by rule $3$. It is easy to show that this $T$ is indeed an $\mathrm{EC}(\mathcal{D})$ and that in every 
leaf $u$ the final triple $p_3p_1p_2$ belongs to ${\cal K}$. Therefore $T$ constructs an equilateral triangle from nothing, by means of only compass and two 
classical arbitrary points.
\eproof

\noindent
The Mohr--Mascheroni theorem says  that compass constructions are as strong as general constructions (see N.~Hungerb\"uhler \cite{hung} and G.\,L.~Alexanderson 
\cite{alex}), thus our above definition of general constructions
is in fact superfluous. It is easy to adapt the proof in \cite{hung} to our model with arbitrary points. The fact that 
$$
|\R|=\mathfrak{c}>|\N|\;,
$$ 
the real numbers have uncountable cardinality, is crucial in the proofs of the Theorems~\ref{impo_unit}, \ref{hilb_thm}, \ref{imp_point}, and \ref{proj_hilb_thm}.

For interest and contrast we illustrate by the next example a~different approach to Euclidean constructions; the example 
is taken from the monograph \cite{schr} of P.~Schreiber. 

\begin{exam}[bisector of a~segment, pp.~140--141 in \cite{schr}]\tec\label{ex_schreiber}
Recall that a~{\em bisector} of two distinct points $p_1$ and $p_2$ is the set of points in the plane $\R^2$ that are equidistant 
to $p_1$ and $p_2$. It is the line perpendicular to the segment $p_1p_2$ and going through its midpoint.
\end{exam}
In this paragraph we translate freely from \cite{schr}. The simple problem, to construct for a~straight segment given by its endpoints by compass and 
straightedge its bisector (``die Mittelsenkrechte''), means in our sense the following. One should give a~constructive proof for the CE (conditioned existential 
proposition)
$$
\mathrm{\forall\,p_1p_2\,\big(p_1\ne p_2\to\exists \,g\,(g\perp L(p_1,p_2)\wedge S(L(p_1,p_2),g)p_1\cong S(L(p_1,p_2),g)p_2)\big)\;,}
$$
by presenting a~uniform flowchart over $(\mathcal{T},\mathcal{K},\mathcal{E})$. Here $\mathcal{T}$ is the plane Euclidean geometry, $\mathcal{K}$ is the 
system of the five CUEs (conditioned univalent existential propositions) corresponding to the operations $\mathrm{L,Z,S_1,S_2,S_3}$, and $\mathcal{E}$,
as it tuns out, can be taken empty. The standard solution of this problem, described as a~uniform flowchart, reads:
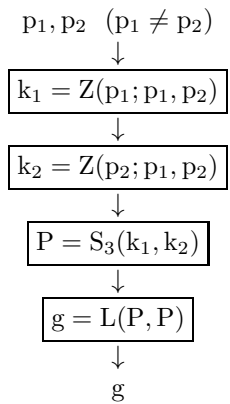
\begin{center}
\begin{tabular}{c}
   $\mathrm{p_1,p_2\ \ (p_1\ne p_2)}$    \\
   $\downarrow$\\
    \fbox{$\mathrm{k_1=Z(p_1;p_1,p_2)}$}  \\
    $\downarrow$\\
    \fbox{$\mathrm{k_2=Z(p_2;p_1,p_2)}$}  \\
    $\downarrow$\\
    \fbox{$\mathrm{P=S_3(k_1,k_2)}$}  \\
    $\downarrow$\\
    \fbox{$\mathrm{g=L(P,P)}$}  \\
    $\downarrow$\\
    $\mathrm{g}$  \\
\end{tabular}
\end{center}
To be proven: under the assumption $p_1\ne p_2$, the circles $k_1,k_2$ are defined, intersect each other, the line going through the intersections 
is perpendicular to $L(p_1,p_2)$ and the intersection of $L(p_1,p_2)$ with this line has the same distances to $p_1$ and $p_2$.

This passage hopefully conveys to the reader some flavor of the approach in \cite{schr}.
\eproof

\noindent
Neither Hilbert's theorem nor \cite{caue} or other literature related to the theorem are mentioned in \cite{schr}. It appears that the problem 
of arbitrary points in Euclidean constructions is outside the scope of interest of \cite{schr}.

So we turn to two formal properties of $\mathrm{EC}(\mathcal{S})$; these results are not very deep but we think that they illustrate nicely 
the above formal definitions.

\begin{prop}[on constructing]\tec\label{prop_constr}
Let ${\cal S}$, ${\cal K}$, and $r$ be as above, and ${\cal K}$ and $r$ be separated. Then an $\mathrm{EC}(\mathcal{S})$ of a~type constructs 
${\cal K}$ from $r$ if and only if an $\mathrm{EC}(\mathcal{S})$ of the same type weakly constructs ${\cal K}$ from $r$. 
\end{prop}
\proof
The implication from constructing to weak constructing is trivial. To show the opposite implication we assume that $T=(r,V,E)$ is a terminating 
$\mathrm{EC}(\mathcal{S})$ of some type that weakly constructs ${\cal K}$ from $r$. We check that for any vertex $u=a_1a_2\ds a_m\in V$ 
and its any letter $a_i\not\in{\cal S}$ that appears also as a~letter in a~word in ${\cal K}$, one can make $v=a_1a_2\ds a_ma_{m+1}$ with $a_{m+1}=a_i$ 
a~child of $u$ according to one of the rules 2--6 in the definition of $\mathrm{EC}(\mathcal{S})$ and
according to the type of $T$ (for $a_i\in{\cal S}$ it is actually also possible). By these one-vertex prolongations we can 
prolong every branch in $T$ so that the resulting $T'$ is a~terminating $\mathrm{EC}(\mathcal{S})$ of the same type as $T$ and constructs ${\cal K}$ from $r$.
Note that $a_i$ appears in $u$ as a result of application of one of the rules 2--6 and not just because it was in $r$ already at the start: ${\cal K}$ and 
$r$ are separated.

If $a_i\in{\cal C}$ then $a_i$ appears in $u$ because $a_i=k(a_j,a_k,a_l)$ for some three points with indices $1\le j,k,l<i$ and with $a_k\ne a_l$ 
and the type of $T$ is not straightedge construction. We can apply on $u$ rule $3$ again, with the same indices $j,k$, and $l$, and make $v$ a child of 
$u$. For $a_i\in{\cal L}$ the argument is the same, now the type of $T$ is not compass construction and we apply again rule $2$. Suppose that $a_i\in{\cal P}$. 
Then $a_i$ appears in $u$ either because it is an intersection point of two distinct curves $a_j$ and $a_k$ with $1\le j,k<i$ or because 
$a_i\in a_{i-1}\in{\cal S}$ is an arbitrary point or because $a_i=k(a_j,a_k,a_k)$ for two points $a_j$ and $a_k$ with $1\le j,k<i$. In each of the three cases 
we can simply repeat $a_i$ as a~degenerated circle by rule 3: $a_{m+1}=k(a_i,a_i,a_i)$. In the first and third case we can alternatively apply again 
the same rule 4 or rule 3, but in the second case the use of a~degenerated circle is unavoidable.
\eproof

Another general transformation of constructions is the following. 
Suppose that ${\cal S}$ and ${\cal S}'$ are nonempty systems of nonempty subsets of $\R^2$ such that ${\cal S}'\sus{\cal S}$ and for every $S\in{\cal S}$ 
there is an $S'\in{\cal S}'$ with $S'\sus S$. Let ${\cal K}\sus({\cal P}\cup{\cal L}\cup{\cal C})^*$ and $r\in({\cal P}\cup{\cal L}\cup{\cal C})^*$, 
not necessarily separated.

\begin{prop}[equivalent models]\tec\label{equi_mode}
In this situation, an $\mathrm{EC}(\mathcal{S})$ of a~type constructs ${\cal K}$ from $r$ if and only if an $\mathrm{EC}(\mathcal{S}')$ of the same type does.
\end{prop}
\proof
If $T'=(r,V,E)$ is an $\mathrm{EC}(\mathcal{S}')$ of some type constructing ${\cal K}$ from $r$, by the assumption it is also an $\mathrm{EC}(\mathcal{S})$ 
of the same type. 

If $T=(r,V,E)$ is an $\mathrm{EC}(\mathcal{S})$ of some type constructing ${\cal K}$ from $r$, we define a rooted tree $T'=(r',V',E')$ with $r'=r$, roughly 
a~rooted subtree of $T$, that is an $\mathrm{EC}(\mathcal{S}')$ of the same type as $T$ and also constructs ${\cal K}$ from $r$. We proceed by induction 
on the {\em height $i\in\N_0$} of a~vertex $u\in V'$, which is the length (the number of edges) of the path from $r'$ to $u$.
Along we also inductively define an injection $$
f\cc V'\to V,\ f(u)=f(a_1a_2\ds a_m)=b_1b_2\ds b_m\;, 
$$
such that $a_i\in{\cal S}'\iff b_i\in{\cal S}$ and 
$a_i\not\in{\cal S}'\Rightarrow a_i=b_i$. For $i=0$ we set $u=r'=r$ and $f(u)=u$. Suppose that $i>0$ and that all vertices in $T'$ with height less than $i$ 
have been already defined, as well as the edges between them and the values of $f$ on them. We consider all vertices $u$ in $T'$ with height $i-1$. If there 
is none, we are done with $T'$ and its definition is at the end. Let $u=a_1a_2\ds a_m\in V'$ have height $i-1$ and be deterministic. If $f(u)$ is a~leaf in $T$, 
we keep $u$ a~leaf in $T'$ too. Else we consider the unique edge $e=(f(u),v)\in E$, $v=b_1b_2\ds b_{m+1}$. If $e$ was not obtained by rule $5$, we add to $V'$ 
the new vertex $v'=a_1a_2\ds a_mb_{m+1}$, to $E'$ the new edge $(u,v')$, and we set $f(v')=v$. If $e$ was obtained by rule $5$ and $b_{m+1}\in{\cal S}$, we take an 
$a_{m+1}\in{\cal S}'$ such that $a_{m+1}\sus b_{m+1}$ and add to $V'$ the new vertex $v'=a_1a_2\ds a_{m+1}$, to $E'$ the new edge $(u,v')$, and we set $f(v')=v$. 
If $u=a_1a_2\ds a_m$ is non-deterministic, then we add to $V'$ the new vertices $v'=a_1a_2\ds a_{m+1}$ for each point $a_{m+1}$ in $a_m\in{\cal S}'$, to $E'$  
the corresponding new edges $(u,v')$, and for each $v'$ we set $f(v')=v$ where $v\in V$ is the unique child of $f(u)$ in $T$ whose last letter is the point 
$a_{m+1}$. The rooted tree $T'$ consists of exactly all vertices and edges obtained when $i$ runs in $\N_0$. Clearly, $T'$ has the same root as $T$ 
and it follows from its inductive definition that it is a~terminating $\mathrm{EC}(\mathcal{S}')$. Since every leaf $u$ in $T'$ coincides, except for the 
letters in ${\cal S}'$, with the leaf $f(u)$ in $T$, it follows that $T'$ constructs ${\cal K}$ from $r$ too. It is clear that $T'$ is of the same type as $T$. 
\eproof

\noindent
The set systems ${\cal S}={\cal O}$ and ${\cal S}'={\cal D}$, of nonempty open sets and of open discs with positive radii, form an example of 
the situation treated by the proposition. In place of $\mathcal{S}'=\mathcal{D}$ we may take any other basis of the Euclidean topology on $\R^2$.
\begin{quote}
    {\em Following \cite{akop_fedo,shen,shen_uspe}, we postulate that the $\mathrm{EC}(\mathcal{D})$, and equivalently the  $\mathrm{EC}(\mathcal{O})$, provide a~rigorous 
model of Euclidean constructions with arbitrary points.}
\end{quote}

The two previous propositions treat the $\mathrm{EC}(\mathcal{S})$ only as data structures. We turn to more substantial results
on $\mathrm{EC}(\mathcal{S})$ and begin with an auxiliary lemma. We say that a set $X\sus\R^2$ is E-{\em closed} if for any ten elements $a,b,c,d,e,f,g,h,i,j$ in $X$ such that $a\ne b$ and $c\ne d$,
each of the sets
$$
l(a,\,b)\cap l(c,\,d),\ l(a,\,b)\cap k(e,\,f,\,g),\ \mathrm{and}\ k(e,\,f,\,g)\cap k(h,\,i,\,j)
$$
that has at most two elements is contained in $X$.

\begin{lem}[on E-closed sets]\label{e_closed}\tec
For every countable set $X\sus\R^2$ there is a~countable set $Y\sus\R^2$ such that $X\sus Y$ and $Y$ is {\rm E}-closed. 
\end{lem}
\proof
We set $Y_0=X$. If $Y_0,Y_1,\ds,Y_n$, $n\in\N_0$, have been already defined, we define $Y_{n+1}$ to be the union of all of the one- and two-element
sets of intersection points displayed above, for all ten-tuples $a,b,\ds,j\in Y_0\cup Y_1\cup\ds\cup Y_n$ with $a\ne b$ and $c\ne d$. Then we set
$$
Y=\bigcup_{n=0}^{\infty}Y_n\;.
$$
It is easy to see that $Y$ has the stated properties.
\eproof

The next theorem is a~``baby version'' of Hilbert's Theorem~\ref{hilb_thm}; it has a~simple proof nicely illustrating the proof method.
\begin{thm}[non-constructibility of unit length]\tec\label{impo_unit}
Every (classical and not necessarily terminating) $\mathrm{EC}(\mathcal{D})$ $T=(r,V,E)$ with $r=\emptyset$ has a~branch $B$ such that $|pq|\ne1$ for every two distinct points $p,q\in u$ in every 
vertex $u\in B$. Thus no general classical $\mathrm{EC}(\mathcal{D})$ constructs 
$$
{\cal K}=\{ab\in{\cal P}^2\;|\;|ab|=1\}
$$
---\,the unit length\,---\,from nothing. 
\end{thm}
\proof
We prove existence of a~set $X\sus\R^2$ with three properties:
\begin{enumerate}
    \item $X$ is dense in $\R^2$, every disc $D\in{\cal D}$ intersects $X$. 
    \item For every two distinct points $p,q\in X$, the distance $|pq|\ne1$.
    \item $X$ is E-closed. 
\end{enumerate}
With such set $X$ it is easy to show that every $\mathrm{EC}(\mathcal{D})$ $T=(\emptyset,V,E)$ has the required branch $B=v_1v_2\ds v_m$ or $B=v_0v_1\ds\;$.
We set the first vertex of $B$ to be  $r=\emptyset$. If the vertex $v_n$, $n\in\N_0$, of $B$ has been already defined, we distinguish the cases of deterministic
and non-deterministic $v_n$. In the former case we finish $B$ with $v_n$ if it is a leaf, 
and else set $v_{n+1}$ to be the child of $v_n$. For non-deterministic $v_n=a_1a_2\ds a_m$ with $a_m\in{\cal D}$ we set $v_{n+1}=a_1a_2\ds a_{m+1}$ for some 
$a_{m+1}\in a_m\cap X$, which is possible by property $1$ of $X$. It is clear by the definition of $B$ and $T$ and by property $3$ of $X$ that for every 
$v\in B$, every letter in $v$ that is a point lies in $X$. Thus by property $2$ of $X$, no vertex in $B$ contains two point letters with distance $1$. 

We have to show that a~set $X$ with properties 1--3 exists. By Lemma~\ref{e_closed}, there is a countable and E-closed set $X'\sus\R^2$ with
$X'\supset\Q\times\Q$. This $X'$ is also dense. Thus $X'$ has properties $1$ and $3$ but not property $2$. To achieve it, we modify $X'$ so that
all distances $1$ between its points are destroyed but properties $1$ and $3$ are preserved. We consider the countable set
$$
M=\{1/|ab|\;|\;a,\,b\in X',\,a\ne b\}\sus\R\;.
$$
As $(0,+\infty)$ is an uncountable set, there is a positive real number $\al$ such that $\al\not\in M$. We define
$$
X:=\al X'=\{(\al x,\,\al y)\;|\;(x,\,y)\in X'\}\;.
$$
We show that $X$ has properties 1--3. For every $D\in{\cal D}$, $\al^{-1}D\in{\cal D}$.
So $X$ is dense because $X'$ is dense\,---\,$X$ has property $1$. For every $a,b,c,d\in\R^2$ with $a\ne b$ we have 
$$
\al^{-1}l(a,\,b)=l(\al^{-1}a,\,\al^{-1}b)\ \mathrm{and}\ \al^{-1}k(a,\,c,\,d)=k(\al^{-1}a,\,\al^{-1}c,\,\al^{-1}d)\;,
$$
which implies that $X$ is E-closed because $X'$ is  E-closed. So $X$ has property~$3$. To check that $X$ has property $2$, we assume for contrary that $|ab|=1$ for 
two points $a,b\in X$. But then the points $a'=\al^{-1}a$ and $b'=\al^{-1}b$ lie in $X'$, $|a'b'|=\al^{-1}$, and $\al=|a'b'|^{-1}\in M$, contrary to the definition of $\al$. 
Thus $X$ has all properties 1--3.
\eproof

\noindent
When we replace classical arbitrary points ${\cal S}={\cal D}$ 
with ${\cal S}={\cal U}$, arbitrary points determined by proper horizontal segments, it  becomes possible to construct unit length by compass and straightedge from nothing.

\begin{prop}[constructibility of unit length]\tec\label{po_unit}
Let
$$
{\cal K}=\{ab\in{\cal P}^2\;|\;|ab|=1\}\ \mathrm{and}\ r=\emptyset\;.
$$
There exists a general {\rm EC}$({\cal U})$ $T$ constructing ${\cal K}$ (a unit length) from nothing.
\end{prop}
\proof
By now we may describe $T$ less formally than in Example~\ref{example}. By selecting four appropriate ${\cal U}$-arbitrary points we construct (by straightedge) the two lines $y=0$ and $y=1$.
Then, by intersecting the line $y=0$ with lines going through appropriate ${\cal U}$-arbitrary points, we select two distinct arbitrary points $a,b\in(y=0)$. Finally,
we construct by the standard construction (recalled in Example~\ref{ex_schreiber}) the bisecting line $\ell$ of the points $a$ and $b$ (for which we need compass). The intersections $\ell\cap(y=0)$ and  $\ell\cap(y=1)$ are two points with distance $1$.
\eproof

Interestingly, with $\mathcal{U}$-arbitrary points it is also possible to give a~construction deemed impossible in Hilbert's theorem.
\begin{prop}[``refutation'' of Hilbert's theorem]\tec\label{hilb_inval}
Let $k\in{\cal C}$ be a circle in the plane, $c\in{\cal P}$ be its center, and 
$$
{\cal K}=\{c\}\ \mathrm{and}\ r=k\;.
$$
There exists a~straightedge $\mathrm{EC}(\mathcal{U})$ $T$ constructing ${\cal K}$ (the center $c$ of $k$) from $r$ (the given circle~$k$).
\end{prop}
\proof
We again describe $T$ informally. By selecting six appropriate ${\cal U}$-arbitrary points we construct by straightedge three horizontal lines $y=t_1$, $y=t_2$, and $y=t_3$ such that $t_i\in\R$, $t_1$ is the $y$-coordinate of $c$, $t_1<t_2<t_3$, and $t_3-t_1$ is smaller than the radius of $k$. We construct the six intersection points $p_1,\ds,q_3$ 
of these lines with the circle $k$ and denote their coordinates as
$$
p_i=(x_i,\,t_i),\ q_i=(x_i',\,t_i),\ x_i<x_i',\ i=1,\,2,\,3\;.
$$
We construct by straightedge the four lines $\ell_i=l(p_1,q_{i+1})$ and $\ell_i'=l(p_{i+1},q_1)$, $i=1,2$. We construct the two intersection points
$a=\ell_1\cap\ell_1'$ and $b=\ell_2\cap\ell_2'$. Finally, we construct by straightedge the line $\ell_3=l(a,b)$. The intersection point
$\ell_3\cap(y=t_1)=c$, the center of $k$. We only used straightedge, never compass. But we used the stronger ${\cal U}$-arbitrary points,
not the classical ones.
\eproof

\noindent
Thus an argument supporting Hilbert's theorem but not specifying precisely (on the level of Definition~\ref{defi_EC}, say) selection of
arbitrary points may at best be only an idea that possibly may (or may not) lead to a~rigorous proof. As far as we know, this is the case with all proofs of Hilbert's theorem
and related results prior to \cite{akop_fedo}. We quote from the literature one such argument and discuss its shortcomings. 
One can choose from at least three sources: R.~Courant and H.~Robbins \cite[p. 152]{cour_robb} and M.~Kac and S.\,M.~Ulam \cite{kac_ulam} and \cite[p. 18]{kac_ulam_dove} in English, 
and H.~Rademacher and O.~Toeplitz \cite[pp. 151--152]{rade_toep} in German. For the first quote see also \cite{shen}. We choose the last 
quote as it is most detailed and, by time of origin and language, perhaps closest to the original argument of D.~Hilbert. We remark that the beginning of the
second quote \cite[p. 18]{kac_ulam_dove}  incorrectly attributes Hilbert's theorem to J.~Steiner. This arose probably by confusion with the 
theorem of J.~Steiner that every Euclidean construction by compass and straightedge can be performed only by straightedge, if one circle together with its center are 
given (\cite{alex,caue}). Hilbert's theorem shows that this center is indispensable. We also remark that this author first learned about Hilbert's theorem 
and its proof in \cite[p. 25]{kac_ulam_cz}, the translation of \cite{kac_ulam} in Czech, in the late 1980s. Then he of course did not notice anything suspicious.
After the quote we translate it to English. In Chapter 21.3 of \cite{rade_toep} one can read on pp. 151--152 the following.
\begin{quote}{\small
3. Nehmen wir an, wir h\"atten zu einem gezeichnet vorliegenden Kreis durch blo\ss e Benutzung des Lineals nach einem gewissen Verfahren den Mittelpunkt konstruieren. 
Man h\"atte also gerade Linien gezogen, die den Kreis oder einander schneiden und h\"atte gewisse Schnittpunkte durch gerade Linien verbunden. Da hierbei ein Punkt nur
fixiert werden kann durch gerade Linien, auf denen er liegt, so w\"are also schlie\ss lich der Mittelpunkt als der Schnittpunkt zweier Geraden in diesem Verfahren 
aufgetreten. Die so erzielte Figur best\"ande also aus dem gegebenen Kreis und einigen geraden Linien, von denen zwei sich im gesuchten Kreismittelpunkt schneiden.

Wir werden nun eine besondere {\em Abbildung} dieser Figur studieren, eine Abbildung, die zun\"achst den Kreis wieder in einen Kreis \"uberf\"uhrt, jede gerade Linie
in eine gerade Linie und jeden Schnittpunkt wieder in den Schnittpunkt der entsprechenden Linien. Solcher Abbildungen gibt es nat\"urlich sehr viele; z. B. w\"are jede
\"ahnliche Vergr\"o\ss erung oder Verklei-\\nerung der Figur eine solche. Aber gerade mit \"ahnlichen Abbildungen ist uns
f\"ur unseren Zweck nich gedient. Wir werden Vielmehr eine solche Abbildung angeben, die zwar unseren Kreis als Kreis und jede Gerade als Gerade erh\"alt, aber doch 
die Figur v\"ollig verzerrt, vor allem den Kreismittelpunkt in einen Bildpunkt \"uberf\"uhrt, der gewi\ss\ nich der Mittelpunkt des Bildkreises ist. 

Wenn wir eine solche Abbildung angeben k\"onnen, sind wir schon fertig mit unserem Beweis. Denn in der Tatt: die Bildfigur mag sich von der Originalfigur noch so sehr
unterscheiden, in bezug auf die als m\"oglich angenommene Konstruktion sind beide Figuren v\"ollig gleichberechtigt. Jeden Schritt der Konstruktion in der 
Originalfigur, etwa das Ziehen einer Geraden, das Aufsuchen eines Schnittpunktes oder das Verbinden zweier Schnittpunkte durch eine Gerade, k\"onnten wir auch, da der
Kreis und jede Gerade und jeder Schnittpunkt sich im Bilde wiederfinden, in derselben Reihenfolge in der Bildfigur ausf\"uhren. Da aber nach Voraussetzung der Mittelpunkt
des Originalkreises {\em nicht} auf den Mittelpunkt des Bildkreises abgebildet ist, so kann die Konstruktion in der Bildfigur nicht zum Ziele gef\"uhrt haben: zu den Geraden,
die sich in der Originalfigur im Mittelpunkt des Kreises schneiden sollten, geh\"oren Bildgeraden, deren Schnittpunkt vom Mittelpunkt des Bildkreises verschieden ist.
Obgleich also auch in der Bildfigur Schritt f\"ur Schritt die angenommene Konstruktionvorschrift erf\"ullt geworden ist, hat sie doch nicht die Auffindung des 
Kreismittelpunktes geleistet. Das ist aber ein Widerspruch gegen den Sinn einer Konstruktionmethode. Also kann es eine solche gar nicht geben: mit dem Lineal allein ist
die Konstruktion des Mittelpunktes eines ohne Mittelpunkt gegebenen Kreises unausf\"uhrbar. 

F\"ur den Fall {\em zweier} Kreise wird unser Beweis nachher ganz analog verlaufen. [emphasizes in the original]}
\end{quote}
Here is our imperfect translation; neither German nor English is mother tongue of this author. Word orders in the three languages are not easy to reconcile. 
\begin{quote}{\small
3. Let us suppose that we have constructed for a~given circle, lying drawn before us, by a~certain procedure and by means of a~mere straightedge, the center. One has drawn straight lines that 
intersect the circle or one another and one has connected specified intersection points by straight lines. As by this a~point can be determined only by the lines on 
which it lies, at the end of the procedure the center has appeared as the intersection point of two lines. Thus obtained picture consists therefore of the given circle
and some straight lines, of which two intersect in the sought-for center of the circle.

We will study a special {\em mapping/transformation} of this picture, a mapping that firstly maps the circle again to a circle, every line to a line and every intersection again to the
intersection of the corresponding lines. Naturally, there are very many of such mappings: for example, such is every magnifying or downsizing similarity of the picture. But exactly 
similarity mappings cannot serve for our goal. We will give even such mapping that on the one hand preserves our circle as a circle and every line as a line, but on the 
other hand completely deforms the picture,  first of all it sends the center of the circle to an image point that surely is not the center of the image circle.

If we can give such mapping, we are done with our proof. Indeed: the image picture cannot substantially differ from the original one, with regard 
to the supposedly possible construction
both pictures have completely equal rights. Every step of the construction in the original picture, like drawing a line, finding an intersection or connecting 
two intersections 
by a line, we could also, since the circle and every line and every intersection appear again in the image, perform in the same order in the image picture. But since by the
assumption the center of the original circle does {\em not} map on the center of the image circle, the construction in the image picture cannot lead to the goal:  to the lines that in the original picture should intersect in the center of the circle,  image lines correspond whose intersection differs from the center of the image circle. Although also
in the image picture the supposed construction recipe has been followed step by step, it has not succeeded in finding the center of the circle. But this is a contradiction 
with the sense 
of a construction method. Thus no such construction can be at all: with straightedge alone the construction of the center of a circle, that is given without the center, is unperformable.

In the case of {\em two} circles our proof will be accomplished later quite similarly. 
}
\end{quote}
We rise three objections to this argument. {\em First objection.} Formulations are vague, the procedure (''das Verfahren'') is not formally defined. 
This is actually the main problem. Could the authors reply to us, they would probably say that almost all arguments in mathematics are informal and that in this 
case it would present no problem to formalize the procedure as a~sequence of precisely defined steps, etc.
{\em Second objection.} A~strange feature of the quoted argument is that it is worded as if the construction were deterministic, arbitrary points 
are not mentioned at all; we return to it in the next section. Similarly,  R.~Courant and H.~Robbins \cite[p. 152]{cour_robb} describe the hypothetical 
construction as if it were deterministic, without mentioning arbitrary points. Only M.~Kac and S.~M.~Ulam \cite[p. 18]{kac_ulam_dove} write: ``For example, 
a~step may call for choosing two arbitrary points on the circumference of the circle and joining them by a straight line.'' {\em Third objection} pinpoints
the error, we think, as an error in intuition: a~logical claim is presented as a~sure thing, without any justification or proof. We mean the sentence ``Denn in der Tatt: 
die Bildfigur mag sich von der Originalfigur noch so sehr unterscheiden, in bezug auf die als m\"oglich angenommene Konstruktion sind beide Figuren v\"ollig gleichberechtigt.''\,---\,``Indeed: the image picture cannot substantially differ from the original one, with regard to the supposedly possible construction 
both pictures have completely equal rights.'' By this the authors mean, probably, that a~mapping/transformation that sends the given circle to a~circle, and any 
line to a~line, and hence sends every intersection point of two lines or of a~line and the circle again to an intersection point of the image objects, 
{\em also preserves the ``sense of a~construction method''} and has to map the intersection of the two lines, obtained in the construction so that they intersect in 
the center of the given circle, to the center of the image circle. But if one thinks about it a~while, especially after the experience of Proposition~\ref{hilb_inval},
one does not see any clear reason why this claim should hold at all. We call this logical gap the {\em transformation fallacy}. 

On a~positive note we have to
say that the deforming transformation discussed in all three quotes does eventually lead to correct proofs. But, as it turns
out, correct proofs need a~whole uncountable family of such transformations (as in the proof of Theorem~\ref{impo_unit}), one does not suffice.

We proceed to a~correct proof of Hilbert's theorem. Correct proofs, for formulations differing from Theorem~\ref{hilb_thm}, were given
already in \cite{akop_fedo} and \cite{shen}. Our proof runs in the framework of $\mathrm{EC}(\mathcal{D})$ and is quite detailed\,---\,given the history 
of fallacious proofs of Hilbert's theorem and related results, one has to be careful.

\begin{thm}[rigorous Hilbert's theorem]\tec\label{hilb_thm}
Let $k\in\mathcal{C}$ be a~circle in the plane, $c\in\mathcal{P}$ be its center, and let $\mathcal{K}=\{c\}$.
Then every straightedge classical (not necessarily terminating) $\mathrm{EC}(\mathcal{D})$ $T=(r,V,E)$ with $r=k$ has a~branch $B$ such that $c\not\in u$ for any vertex $u\in B$. Thus no straightedge classical $\mathrm{EC}(\mathcal{D})$ 
constructs $\mathcal{K}$ (the center $c$ of $k$) from $r$ (the given circle $c$).
\end{thm}
The proof method is the same as for Theorem~\ref{impo_unit} and is based on a~set $X\sus\R^2$ with the following three properties.
\begin{enumerate}
    \item $X$ is dense in $\R^2$. 
    \item $c\not\in X$.
    \item $X$ is {\em {\rm H}-closed}, for every quadruple $a,b,d,e\in X$ of non-colinear points with $a\ne b$ and $d\ne e$ one has that
    $$
    l(a,\,b)\cap k\sus X\ \mathrm{and}\ l(a,\,b)\cap l(d,\,e)\sus X\;.
    $$
\end{enumerate}
Assuming that such set $X$ exists, we proceed as in the proof of Theorem~\ref{impo_unit} and inductively define in any given straightedge $\mathrm{EC}(\mathcal{D})$ a~branch 
$B$ such that in every vertex $u\in B$ every point letter lies in $X$ and thus differs from $c$. It remains to prove that $X$ exists.

We define $X$ by means of the above mentioned transformation/map preserving the set of lines and the given circle but not its center. A~technical complication is 
that the maps we will use now, unlike the maps $x\mapsto\al x$ in the proof of Theorem~\ref{impo_unit}, are not defined everywhere and are not onto. 
We work in the affine plane $\R^2$ and (unlike in \cite{akop_fedo,shen}) do not employ any projective elements. A~projective version of Hilbert's 
theorem is presented in Theorem~\ref{proj_hilb_thm}. The required maps are given in Proposition~\ref{map_fp} and after proving it we conclude the 
proof of  Theorem~\ref{hilb_thm}. 

For two different lines $\ell,\ell'\in\mathcal{L}$ we write $\ell\setminus\ell'$ for the {\em deleted line $\ell$}, the line $\ell$ with the possible intersection 
point $\ell\cap\ell'$ deleted. For $\ell\in{\cal L}$ we set ${\cal L}_{\ell}=\{\ell'\setminus\ell\;|\;\ell'\in{\cal L},\ell'\ne\ell\}$. For 
$\kappa=\ell\setminus\{p\}$, where $\ell\in{\cal L}$ and  $p\in{\cal P}$, we set $i(\kappa)=\emptyset$ if $\kappa=\ell$ and $i(\kappa)=p$ else. 
In the following proposition we follow Gy.~Strommer \cite[p. 97]{stro}. 

\begin{prop}[Strommer's map]\tec\label{stro_map}
Let $\ell_0=(x=0)$ be the $y$-axis, $a>1$ be a~real number, $k$ be the circle $$(x-a)^2+y^2=a^2-1$$ with the center $c=(a,0)$, 
and let $f$ be the map  
$$
f\cc\R^2\setminus\ell_0\to\R^2\setminus\ell_0,\ f(x,\,y)=(x',\,y'):=(1/x,\,y/x)\;.
$$
Then the following hold. 
\begin{enumerate}
\item The map $f$ is an involution ($f=f^{-1}$) and a~homeomorphism (a~bijection continuous in both directions).
\item One has that $\ell_0\in{\cal L}$ and $\ell_0\cap k=\emptyset$.
\item For every line $\ell\in{\cal L}\setminus\{\ell_0\}$ one has that $f(\ell\setminus\ell_0)\in{\cal L}_{\ell_0}$.
\item Let $\ell,\ell'\in{\cal L}\setminus\{\ell_0\}$ be distinct lines. Then 
$$
\text{$i(f(\ell\setminus\ell_0))=\emptyset\iff\ell\parallel\ell_0$, 
and $\ell\parallel\ell'\iff i(f(\ell\setminus\ell_0))=i(f(\ell'\setminus\ell_0))$}\;.
$$
\item It is true that $f(k)=k$ and $f(c)\ne c$.
\end{enumerate}
\end{prop}
\proof
Properties 1 and 2 are easy to check. Since
$$
f((\al x+\be y+\ga=0)\setminus\ell_0)=(\al+\be y'+\ga x'=0)\setminus\ell_0\;,
$$
property~3 follows. The map $f=f^{-1}$ transforms two distinct vertical parallel lines $x+\ga=0$ and $x+\ga'$ different from $\ell_0$, 
$\ga\ne\ga'$ and $\ga,\ga'\ne0$, in two distinct vertical parallels $x'+1/\ga=0$ and $x'+1/\ga'$ different from $\ell_0$. It transforms two distinct non-vertical 
parallel deleted lines 
$$
(\al x+y+\ga=0)\setminus\ell_0\;\text{and}\;(\al x+y+\ga'=0)\setminus\ell_0,\ \ga\ne\ga'\;, 
$$
in two distinct deleted lines 
$$
(\al+y'+\ga x'=0)\setminus\ell_0\;\text{ and }\;(\al+y'+\ga'x'=0)\setminus\ell_0
$$
``intersecting''
in the point $(0,-\al)\in\ell_0$, and vice versa. Thus property 4 holds. For property 5 we check that the equation for $k$ transforms under $f$ to itself: 
$$
(1/x'-a)^2+(y'/x')^2=a^2-1
$$ 
is equivalent to $(1-ax')^2+(y')^2=(ax')^2-(x')^2$ which is indeed equivalent to $(x'-a)^2+(y')^2=a^2-1$. Also, $f(c)=f((a,0))=(1/a,0)\ne(a,0)=c$.
\eproof

\noindent
Strommer's algebraic definition of $f$ is more straightforward and easier to work with than the definition by the central projection from 
a~point between two planes in $\R^3$, as given in \cite[Chapter 21.4]{rade_toep} and elsewhere. In the next proposition we generate by the map $f$ 
the required uncountable family of deforming transformations

\begin{prop}[required maps]\tec\label{map_fp}
Let $\ell_0$, $k$, and $c$ be as in the previous proposition. There exists a circle $k_0$, concentric with $k$, such that for every point 
$p\in k_0$ there is a homeomorphism 
$$
f_p\cc\R^2\setminus\ell_0\to\R^2\setminus\ell_p
$$
with the following properties. 
\begin{enumerate}
\item One has that $\ell_p\in{\cal L}$ and $\ell_p\cap k=\emptyset$.
\item For every line $\ell\in{\cal L}\setminus\{\ell_p\}$ one has that $f_p^{-1}(\ell\setminus\ell_p)\in{\cal L}_{\ell_0}$, and similarly for the map $f_p$.
\item Let $\ell,\ell'\in{\cal L}\setminus\{\ell_p\}$ be distinct lines. Then
$$
\text{$i(f_p^{-1}(\ell\setminus\ell_p))=\emptyset\iff\ell\parallel\ell_p$} 
$$
and 
$$
\text{$i(f_p^{-1}(\ell\setminus\ell_p))=i(f_p^{-1}(\ell'\setminus\ell_p))\iff\ell\parallel\ell'$}\;.
$$
\item It is true that $f_p^{-1}(k)=k$ and $f_p(c)=p$.
\end{enumerate}
\end{prop}
\proof
Let $f$ be the map in Proposition~\ref{stro_map}. We set $k_0$ to be the circle with the center $c\in\mathcal{P}$ and radius $|c\,f(c)|$. For $p\in k_0$ we
define $f_p=\varphi\circ f=\varphi(f)$, where $\varphi$ is the rotation of $\R^2$ around $c$ moving $f(c)$ to $p$, and $\ell_p=\varphi(\ell_0)$. Using 
that $f_p^{-1}=f^{-1}\circ\varphi^{-1}=f
\circ\varphi^{-1}$ and the properties of $f$ in Proposition~\ref{stro_map}, it is not hard to check the properties 1--4 in the present proposition.
\eproof

\noindent
{\em Proof of Theorem~\ref{hilb_thm}. }First we assume that $k$ is as in the two previous propositions, say for $a=2$ ($a$ has to be algebraic), and at the end 
we  explain extension to any circle $k'$. To define the set $X$ with the properties 1--3, we first set $X'\sus\R^2$ to be the set of points with algebraic coordinates: $(\al,\be)\in X'$ iff $\al$ and $\be$ are roots of nonzero polynomials with rational coefficients. We use algebraic numbers because we need that
the coordinates of points in $X'$ be closed to addition, alternatively we could define $X'$ as in the proof of Theorem~\ref{impo_unit} by some variant of Lemma~\ref{e_closed}.
It is clear that $X'$ is countable and 
has properties~$1$ and $3$, but not property $2$. To achieve property~2 and at keep properties~$1$ and $3$, we transform $X'$ by a~map $f_p$ 
from Proposition~\ref{map_fp}, for a~point $p\in k_0\setminus X'$. Such $p$ exists because $|k_0|=\mathfrak{c}>|X'|=\aleph_0$.
We set
$$
X=f_p^{-1}(X'\setminus\ell_p)\cup\{\ell_0\cap l(p_1,\,p_2)\;|\;p_1,\,p_2\in f_p^{-1}(X'\setminus\ell_p),\,p_1\ne p_2\}
$$
and denote by $Y\sus\ell_0$ the second set in the union. This a~bit complicated definition of $X$ reflects the facts that the maps $f_p$ are not everywhere defined and transform pairs of ``intersecting'' deleted lines in pairs of parallel deleted lines.

The set $X$ is dense in $\R^2$ (property $1$) because $X'$ is dense, $f_p$ is a~homeomorphism and $\ell_0$ and $\ell_p$ are nowhere dense in $\R^2$.
If $c\in X$ then (since $c\not\in Y$) by property~$4$ in Proposition~\ref{map_fp} we would have $p=f_p(c)\in X'$, contrary to the selection of $p$. Thus $X$ has property~$2$.  
We check in detail that $X$ has property~$3$. Let $p_i\in X$ for $i=1,\ds,4$ be four noncolinear points, $p_1\ne p_2$ and $p_3\ne p_4$, and let $\ell_i=l(p_{2i-1},p_{2i})$ for $i=1,2$. First we check that $\ell_1\cap k\sus X$. 
If $p_1,p_2\in X\setminus Y$ then the line
$$
\kappa=l(f_p(p_1),\,f_p(p_2))
$$
goes through two distinct points in $X'\setminus\ell_p$ and $f_p^{-1}(\kappa\setminus\ell_p)=\ell_1\setminus\ell_0$. By property~$3$ of 
$X'$, $\kappa\cap k\sus X'\setminus\ell_p$. Thus by properties~2 and 4 in Proposition~\ref{map_fp}, 
$$
\ell_1\cap k=f_p^{-1}(\kappa\cap k)\sus f_p^{-1}(X'\setminus\ell_p)\sus X\;.
$$ 
If $p_1\in Y$ and $p_2\in X\setminus Y$ then by the definition of $Y$, $p_1\in l(p_5,p_6)$ for two distinct points $p_5,p_6\in X\setminus Y$. 
If $l(p_5,p_6)=\ell_1$, we are in the previous case. Else we consider the distinct lines $\lambda$, where $\lambda\setminus\ell_p=f_p(\ell_1\setminus\ell_0)$,  
and $\kappa=l(f_p(p_5),f_p(p_6))$. By property 3 in Proposition~\ref{map_fp}, $\lambda\parallel\kappa$. Besides $f_p(p_5)$ and $f_p(p_6)$ there 
are on $\kappa$ infinitely many other points from $X'$ and we can take two distinct of them, $q_1$ and $q_2$, such that
$$
q_3=f_p(p_2)+q_2-q_1\in\lambda\cap(X'\setminus\ell_p)
$$
---we use that algebraic numbers are closed to addition and subtraction. Thus $\ell_1$, as $\ell_1\setminus\ell_0=f_p^{-1}(\lambda\setminus\ell_p)$, goes through two 
distinct points $p_2$ and $f_p^{-1}(q_3)$ in $X\setminus Y$ and we are again in the previously discussed case. If $p_1,p_2\in Y$ then $\ell_1=\ell_0$ and 
$\ell_1\cap k=\emptyset$ by property 2 in Proposition~\ref{stro_map}. 

We check that $\ell_1\cap\ell_2\sus X$. We assume that $\ell_1\cap\ell_2=p'\in{\cal P}$ and show that $p'\in X$. Suppose that $p'\in\ell_0$. 
Not all of $p_1,\ds,p_4$ lie in $\ell_0$, say $p_1\in X\setminus Y$. If $p_2=p'$ then $p'\in X$. If $p_2\ne p'$ then also $p_2\in X\setminus Y$,
and  by the definition of $Y$ we have $p'\in Y$ and $p'\in X$. Suppose that $p'\not\in\ell_0$. Then $p_1$ or 
$p_2$ is in $X\setminus Y$ and $p_3$ or $p_4$ is in $X\setminus Y$, and we may suppose that it holds for $p_1$ and $p_3$. If $p_2\in Y$, we deduce as before 
(by adding to $f_p(p_1)$ an algebraic vector) that $\ell_1\setminus\ell_0=f_p^{-1}(\kappa_1\setminus\ell_p)$ for a line $\kappa_1$ going through two distinct points in $X'\setminus\ell_p$. 
If $p_2\in X\setminus Y$, it holds too and $\kappa_1$ goes through the points $f_p(p_1)$ and $f_p(p_2)$. We define the line $\kappa_2$ for $\ell_2$ 
in the analogous way. Then
$$
p'=f_p^{-1}(\kappa_1\cap\kappa_2)\in f_p^{-1}(X'\setminus\ell_p)\sus X
$$
because $X'$ has property $3$. Thus $X$ has all three properties $1$--$3$ and the proof of Theorem~\ref{hilb_thm} is complete, in the case when $k$ 
is as in the previous two propositions with $a=2$. 

If $k'$ is any circle with center $c'$, we consider the the map 
$$
g=s_1\circ s_2\cc\R^2\to\R^2
$$ 
transforming $k$ to $k'$ and sending $c$ to $c'$, where $s_2$ is the shift $p\mapsto p+c'-c$ of $\R^2$ moving $c$ to $c'$ and $s_1$
is the similarity of $\R^2$ centered at $c'$ that sends the radius $r$ of $k$ to the radius $r'$ of $k'$,
$$
s_1(p)=c'+\frac{r'}{r}(p-c')\;.
$$
It is easy to see that the set $g(X)\sus\R^2$ has with respect to the circle $k'$ properties $1$--$3$. We are therefore done also in the general case. 
\eproof

The simplest nontrivial construction problem is to obtain from nothing by allowed means one prescribed point, which may be taken to be the origin. 

\begin{thm}[non-constructible point]\tec\label{imp_point}
Every classical (not necessarily terminating) $\mathrm{EC}(\mathcal{D})$ $T=(\emptyset,V,E)$ has a~branch $B$ such that $(0,0)\not\in u$ for every vertex 
$u\in B$. Thus no general classical $\mathrm{EC}(\mathcal{D})$ constructs
$$
{\cal K}=\{(0,\,0)\}
$$
(the origin) from nothing. 
\end{thm}

\noindent
The proof is very similar to that of Theorem~\ref{impo_unit} and we omit it. Note that one cannot use exactly the same transformations $x\mapsto\al x$ because 
$(0,0)$ is their fixed point. Instead one can use transformations $p\mapsto\al(p-(0,1))$ or $p\mapsto p+a$, for real $\al>0$ and $a\in\R^2$.

This author believed for some time that the previous theorem also holds with ${\cal U}$ in place of ${\cal D}$. Futile attempts to prove it led eventually
to the opposite conclusion presented in the next theorem. In contrast to Propositions~\ref{po_unit} and \ref{hilb_inval}, the construction is now non-uniform. 

\begin{thm}[constructible point]\tec\label{po_point}
Let ${\cal K}=\{(0,0)\}$, $r=\emptyset$, and ${\cal U}$ be the set system of proper horizontal segments in the plane. There exists a general $\mathrm{EC(\mathcal{U})}$ 
 $T$ constructing ${\cal K}$ (the origin) from nothing. 
\end{thm}
\proof
Since $(0,0)\in(y=0)$, it suffices to construct another line containing the origin. The main idea of the construction is this: If $q$ is a~point outside a~line $\ell$, then for {\em every} two 
distinct points $p,p'\in\ell$,
$$
k(p,\,p,\,q)\cap k(p',\,p',\,q)=\{q,\,\overline{q}\}\;,
$$
where $\overline{q}$ is the mirror image of $q$ with respect to the line $\ell$. 

Now we describe the mechanism of the construction and give an informal description of $T$ in the next paragraph.
Let 
$$
q_i=(q_{i,x},\,i)\in\mathcal{P},\ i=1,\,2\;, 
$$
be two points whose $x$-coordinates satisfy $2q_{1,x}>q_{2,x}>q_{1,x}>0$. The line $\kappa=l(q_2,q_1)$ then intersects 
the $y$-axis in a~point with negative $y$-coordinate. Let $\al\in(0,\frac{\pi}{2})$ be the angle at the vertices $q_i$, determined by the right-going semi-lines 
of $\kappa$ and of the line $y=i$. Let $\ell$ be the line going through $(0,0)$ that subtends in the fourth quadrant with the positive semi-axis $x$ the angle $\al$, 
and consider the point 
$$
b=\ell\cap\kappa\in\mathcal{P}\;.
$$ 
Finally, let $\ell'$ be the line $y=b_y$ where $b_y<0$ is the $y$-coordinate of $b$.
It follows that both angles at the vertex $b$ determined by the right-going semi-lines of $\kappa$, of $\ell'$, and of $\ell$ are equal to $\al$. 
Thus the mirror images $\overline{q_i}$ of $q_i$, $i=1,2$, with respect to the line $\ell'$ lie on the line $\ell$. We can therefore construct $\ell$ as the line determined
by the two points $\overline{q_i}$, where each $\overline{q_i}$ is in turn constructed as above as the other intersection of two circles going through $q_i$ and 
with distinct centers on $\ell'$.

By now the description of $T$ should be clear. First we construct by means of ${\cal U}$-arbitrary points the line $y=0$. Then we construct by means of 
${\cal U}$-arbitrary points two points $q_1$ and $q_2$ whose coordinates satisfy the above conditions. We construct by means of ${\cal U}$-arbitrary points the line $\ell'=(y=b_y)$ where the 
point $b$ depends on $q_1$ and $q_2$ and is defined above\,---\,this is the non-uniform part of the construction because we do not have 
complete control over the position of $b$. We construct by means of ${\cal U}$-arbitrary points two distinct points $p,p'\in\ell'$. We construct the 
other intersection $\overline{q_i}$ of the pair of circles going through $q_i$ and with centers $p$ and $p'$. Finally, we draw the line 
$\ell=l(\overline{q_1},\overline{q_2})$ and get the origin as the intersection 
$$
(0,\,0)=(y=0)\cap\ell\;.
$$
\eproof

Many questions on $\mathrm{EC}(\mathcal{S})$ offer themselves. For instance, what if $\mathcal{S}$ consists of measurable subsets of $\R^2$ with positive measure?
What if $\mathcal{S}$ consists of the ``crosses''
$$
\{((a-\ep,\,a+\ep)\times\{b\})\cup(\{a\}\times(b-\ep,\,b+\ep))\;|\;a,\,b,\,\ep\in\R,\,\ep>0\}?
$$
We hope to tackle these and related questions in \cite{klaz}.

\section{Three variations on Hilbert's theorem}

Unlike the proof of the theorem in \cite{kac_ulam}, the above quoted proof in \cite{rade_toep} and the proof in \cite{cour_robb} do not mention arbitrary points. Could they possibly refer to a~deterministic version of the problem? One can cast Hilbert's theorem deterministicly by taking 
as the starting configuration the given circle and several points on it and allowing only deterministic construction steps. In our terminology the deterministic 
version reads as follows.

\begin{thm}[deterministic Hilbert's theorem I]\tec\label{hilb_thm_determ}
Let $k\in{\cal C}$ be a~circle with the center $c\in{\cal P}$ and let ${\cal K}=\{c\}$. Then for every $n\in\N$ there exist $n$ distinct points 
$p_i\in k$, $i=1,2,\ds,n$, such that in every straightedge deterministic (not necessarily terminating) $\mathrm{EC}(\emptyset)$ $T=(r,V,E)$ with 
$r=kp_1p_2\ds p_n$ we have that $c\not\in u$ for every vertex $u\in V$. Thus no straightedge deterministic $\mathrm{EC}(\emptyset)$ constructs 
${\cal K}$ (the center $c$ of $k$) from $r$ (the given circle $k$ plus the $n$ points $p_i$ on it).
\end{thm}

\noindent
Again, the transformation fallacy is clear: the argument (for impossibility of construction of the center) should supposedly work also for this deterministic version of Hilbert's theorem,   
but for many configurations of the points $p_i$ on $k$ a~deterministic straightedge construction of the center $c$ of $k$ of course exists. 
The simplest of them has four points $p_1,\ds,p_4\in k$ such that the two lines $\ell_i=l(p_{2i-1},p_{2i})$, $i=1,2$, cut $k$ in two distinct diameters, then 
$\ell_1\cap\ell_2=c$. We offer a~strengthening of the previous theorem. Recall that if $Y\sus k$ for a~circle $k$, then $Y$ is {\em dense in $k$} 
if $Y\cap D\ne\emptyset$ for any disc $D\in\mathcal{D}$ intersecting $k$.

\begin{thm}[deterministic Hilbert's theorem II]\tec\label{hilb_thm_det_strong}
Let $k\in{\cal C}$ be a~circle with the center $c\in{\cal P}$ and let ${\cal K}=\{c\}$. There exists a~countable set 
$$
Y\sus k
$$ 
that is dense in $k$ and such that for every finite tuple of points $p_i\in Y$, $i=1,2,\ds,n$, in every straightedge 
 deterministic (not necessarily terminating) $\mathrm{EC}(\emptyset)$ $T=(r,V,E)$ with $r=kp_1p_2\ds p_n$ one has that $c\not\in u$ for every vertex $u\in V$. Thus no straightedge 
deterministic $\mathrm{EC}(\emptyset)$ constructs ${\cal K}$ (the center $c$ of $k$) from $r$ (the given circle $k$ plus some $n$ points from $Y$ on it).
\end{thm}
\proof
Recall that the countable set $g(X)$ defined at the end of the proof of Theorem~\ref{hilb_thm} is dense in $\R^2$, does not contain $c$, and is H-closed. 
It follows that $Y=k\cap g(X)$ has the stated properties. 
\eproof

\noindent
We can state this result more strongly if we modify the above definition of a~$\mathrm{EC}(\mathcal{S})$ 
$T=(r,V,E)$ so that the root is an infinite word with length $\omega$ and the other vertices in $V$ are infinite words with lengths $\omega+i$  for some $i\in\N$. The alphabet is as before, points, lines, circles, and elements of $\mathcal{S}$.
The six rules for children of a~vertex and other notions pertaining to $\mathrm{EC}(\mathcal{S})$ are correspondingly modified. Thus in the next theorem and in Theorem~\ref{concr_hilb_thm} the root $r=a_0a_1\ds$ 
in $T$ has length $\omega$. Recall the notation introduced at the beginning: a~child of $r$ is $a_0a_1\ds;a_{\omega}$ and a~child of a~vertex
$$
\text{$a_0a_1\ds;\,a_{\omega}a_{\omega+1}\ds a_{\omega+i}$ with $i\in\N_0$}
$$
is $a_0a_1\ds;a_{\omega}a_{\omega+1}\ds a_{\omega+i+1}$.
We denote the modified Euclidean constructions with $\mathcal{S}$-arbitrary points and infinite configurations by 
$$
\mathrm{EC}_{\infty}(\mathcal{S})\;.
$$

\begin{thm}[deterministic Hilbert's theorem III]\tec\label{hilb_thm_det_strong3}
Let $k\in{\cal C}$ be a~circle with the center $c\in{\cal P}$ and let ${\cal K}=\{c\}$. There exists a~countable set 
$$
Y=\{p_1,\,p_2,\,\ds\}\sus k
$$ 
that is dense in $k$ and such that in every straightedge deterministic (not necessarily terminating) $\mathrm{EC}_{\infty}(\emptyset)$ 
$T=(r,V,E)$ with $r=kp_1p_2\ds$ one has that $c\not\in u$ for every vertex $u\in V$. Thus no straightedge deterministic 
$\mathrm{EC}_{\infty}(\emptyset)$ constructs ${\cal K}$ (the center $c$ of $k$) from $r$ (the given circle $k$ plus the infinitely many points $Y$ on it).
\end{thm}
\proof
The set $Y=k\cap g(X)$ in the previous proof has the stated properties. 
\eproof

\noindent
The result is formally stronger than Theorem~\ref{hilb_thm_det_strong} because it allows constructions using 
all points in $Y$.

With the help of the standard description of rational points on a~circle we give a~relatively explicit example of a~set $Y\sus k$ in the previous theorem. 

\begin{thm}[concrete deterministic Hilbert's theorem]\tec\label{concr_hilb_thm}
Let $k\sus\R^2$ be the circle
$${\textstyle
(x-\frac{3}{2})^2+y^2=\frac{5}{4}\;\text{ with the center $c=(\frac{3}{2},\,0)$}
}
$$
and $Y=\{p_1,p_2,\ds\}\sus k$ be the countable set of points on $k$, given by

\begin{eqnarray*}
Y&=&\left\{p(\al):=\left(\frac{1}{(\be-\frac{3}{2})s'+\ga s+\frac{3}{2}},\,\frac{(\frac{3}{2}-\be)s+\ga s'}{(\be-\frac{3}{2})s'+\ga s+\frac{3}{2}}\right)\;\bigg|\;\al\in\Q\right\}\;,\\
&&\text{where }\be=\frac{2+2\al^2}{5-2\al+\al^2},\ \ga=\frac{1+4\al-\al^2}{5-2\al+\al^2},\ s=\sin1\;\text{ and }\;s'=\cos1\;.
\end{eqnarray*}
The set $Y$ is dense in $k$ and in every straightedge deterministic (not necessarily terminating) $\mathrm{EC}_{\infty}(\emptyset)$ 
$T=(r,V,E)$ with $r=kp_1p_2\ds$ we have that  $c\not\in u$ for every vertex $u\in V$. Thus no straightedge deterministic 
$\mathrm{EC}_{\infty}(\emptyset)$ constructs ${\cal K}=\{c\}$ (the center $c$ of $k$) from $r$ (the given circle $k$ plus the infinitely many points $Y$ on it).
\end{thm}
\proof
We take the circle $k$ in Proposition~\ref{stro_map} for $a=\frac{3}{2}$. The set $Y_0:=k\cap\Q^2$ is then countable 
and dense in $k$. Namely, $Y_0$ consists of the point $(2,1)$ and the other intersections of $k$ with the lines $\ell_{\al}$ going
through the points $(2,1)$ and $(0,\al)$, for $\al$ running in $\Q$. The line $\ell_{\al}$ is determined by the equation $y=x(1-\al)/2+\al$. We find the two
solutions of the system of the previous linear equation and the equation $(x-\frac{3}{2})^2+y^2=\frac{5}{4}$ and get that
$$
Y_0=\left\{\left(\frac{2+2\al^2}{5-2\al+\al^2},\,\frac{1+4\al-\al^2}{5-2\al+\al^2}\right)\;\bigg|\;\al\in\Q\right\}
$$
(the point $(2,1)$ is the double solution for $\al=2$ when $\ell_2$ is tangent to $k$). By an argument as in the proof of Lemma~\ref{e_closed}
we obtain a countable $\mathrm{H}$-closed set $X'\sus\R^2$ such that $X'\supset Y_0$. It follows that the coordinates of the points in $X'$ are algebraic numbers. Then we  
proceed as in the proof of Theorem~\ref{hilb_thm} and get the desired set $Y$ as
$$
Y=f_p^{-1}(Y_0)
$$
where $f_p^{-1}=f^{-1}(\varphi^{-1})=f(\varphi^{-1})$, $f$ is Strommer's map of Proposition~\ref{stro_map} and $\varphi$ is a~rotation of $\R^2$ around 
$c=(\frac{3}{2},0)$ that moves the point $f(c)=(\frac{1}{3/2},\frac{0}{3/2})=(\frac{2}{3},0)$ to a~point $p$ with transcendental coordinates that is surely 
outside $X'$. Counter-clockwise rotation around $c$ by an angle $\theta\in(0,\frac{\pi}{2})$ moves the point $(\frac{2}{3},0)$ to the point
$${\textstyle
p=(\frac{3}{2}-\frac{5}{6}\cos\theta,\,-\frac{5}{6}\sin\theta)\;.
}
$$
Since $\sin\theta$ and $\cos\theta$ are transcendental numbers for any real algebraic number $\theta\ne0$ (A.~Baker \cite[p.~6]{bake}), we may set $\theta=1$. For a~point 
$(x,y)\in k$ we then have
\begin{eqnarray*}
\varphi^{-1}(x,\,y)&=&((x-3/2)s'+ys+3/2,\,(3/2-x)s+ys')\in k\;\text{ and }\\
f(x,\,y)&=&(1/x,\,y/x)\in k\;\text{ where $s=\sin1$ and $s'=\cos1$}\;.
\end{eqnarray*}
For the set $Y=f_p^{-1}(Y_0)=f(\varphi^{-1}(Y_0))$ we therefore get the above displayed description in the statement of the theorem. Considering the set 
$X$ defined from $X'$ as in the proof of Theorem~\ref{hilb_thm} we see that $Y$ has the stated property.
\eproof

\noindent
Here is a~sample of three points $p(\al)\in Y$:
$$
p(-7)=(1.83944\ds,\,-1.06525\ds),\ p(0)=(0.93113\ds,\,0.96249\ds)
$$
and
$$
p(100)=(1.033100\ds,\,-1.01587\ds)\;.
$$
Using an appropriate shift and similarity, we can move the circle $k$ with its center $c$ to any given circle $k'$ and its center $c'$, respectively, 
and get by this concrete deterministic Hilbert's theorem for $k'$ and $c'$. 

It is straightforward to modify the definition of $\mathrm{EC}(\emptyset)$ for uncountable configurations, and to see that if the initial configuration
is a~circle $k$ plus any arc $A\sus k$ with positive length, then one can deterministicly construct the center of $k$ only by straightedge. 
\begin{quote}
{\em If $k$ is a~circle, is there an uncountable set $A\sus k$ (which is possibly dense in $k$) such that one cannot deterministicly construct the center of $k$ only by straightedge, starting from the given circle $k$ plus the points in $A$ on it?}    
\end{quote}

Our third and last variation on Hilbert's theorem is projective. We review points, lines, circles and discs in the real projective 
plane. Then we adapt $\mathrm{EC}(\mathcal{S})$ to projective geometry and
in Theorem~\ref{proj_hilb_thm} state and prove projective Hilbert's theorem.

We work with the following model  $\mathbb{P}_2$ of the {\em real projective plane}:
$$
\mathbb{P}_2=\{p=\{\overline{x},\,-\overline{x}\}\;|\;\overline{x}\in\mathbb{S}\}\;\text{ where }\;\mathbb{S}=\{(x,\,y,\,z)\in\R^3\;|\;x^2+y^2+z^2=1\}\;.
$$
A~{\em projective point $p$} is thus an {\em unordered} pair of antipodal (symmetric one to another with respect to the origin) points on the unit sphere $\mathbb{S}$ which is situated in the ambient Euclidean space $\R^3$. 
For real triples $\overline{x}=(x_1,x_2,x_3)$ and $\overline{y}=(y_1,y_2,y_3)$ we consider the {\em scalar product}
$$
\langle\overline{x},\,\overline{y}\rangle=x_1y_1+x_2y_2+x_3y_3\;.
$$
A~{\em projective line $\ell$} is determined by a~triple $\overline{y}\in\mathbb{S}$ as
$$
\ell=\{\{\overline{x},-\overline{x}\}\in\mathbb{P}_2\;|\;\langle\overline{x},\,\overline{y}\rangle=0\}\;.
$$
It is, roughly, a~main circle on $\mathbb{S}$, an intersection of the unit sphere with the plane in $\R^3$ going through the origin and with the unit normal 
vector $\overline{y}$. A~{\em projective circle $k$} is determined by a~pair $(\overline{y},a)\in\mathbb{S}\times(0,1)$ as 
$$
k=\{\{\overline{x},\,-\overline{x}\}\in\mathbb{P}_2\;|\;\langle\overline{x},\,\overline{y}\rangle=\pm a\}\;.
$$
It is, roughly, the union of a~spherical circle with a~positive radius and its antipode, and is obtained, roughly, as an intersection of $\mathbb{S}$
with a~pair of distinct antipodal parallel planes with distance less than $1$ from the origin. It is also, roughly, an intersection of $\mathbb{S}$ with a~conic 
surface with the vertex in the origin. A~{\em projective (open) disc $D$} is given by
$$
D=\{\{\overline{x},\,-\overline{x}\}\in\mathbb{P}_2\;|\;\pm\langle\overline{x},\,\overline{y}\rangle>a\}\;,
$$
where $\overline{y}$ and $a$ are as above. It is, roughly, an open spherical cap and its antipode. 
We denote the set of all projective discs by $\mathcal{D}_{\mathrm{pr}}$.
The {\em center $c$} of the above projective circle $k$ (and 
of the corresponding projective disc $D$) is the projective point
$$
c=\{\overline{y},\,-\overline{y}\}
$$
with the same $\overline{y}$ as in the definition of $k$. The {\em radius of $k$} is half of the length of the shortest of the four arcs any projective line $\ell$
going through $c$ is divided into by its four intersections with $k$, when we view $k$, $\ell$ and $c$ as subsets/elements of $\mathbb{S}$. For $p,q,r\in\mathbb{P}_2$
we denote by $k(p,q,r)$ the projective circle or the projective point obtained by taking the projective circle with center $p$ and radius equal to the spherical 
distance (defined in the obvious sense) of the projective points $q$ and $r$ (we again treat these objects as subsets/elements of $\mathbb{S}$); for $q=r$ we set 
$k(p,q,r)=p$. Every two distinct projective points are elements of a~unique projective line and, the nice property of projective geometry, every two distinct 
projective lines intersect in a~unique projective point. 

Another model of the real projective plane is the disjoint union
$$
\mathbb{P}_2'=R\cup L\cup P
$$
where $R\sus\R^3$ is the plane of points $(x,y,1)$, $L\sus\R^3$ is the line of points $(x,1,0)$, and  $P=\{(1,0,0)\}$. The bijection 
$$
F\cc\mathbb{P}_2\to\mathbb{P}_2'
$$ 
sends $p=\{\overline{x},-\overline{x}\}$ to the intersection of $\mathbb{P}_2'$ with the ordinary line $l(\overline{x},-\overline{x})$ in $\R^3$ determined
by the ordinary points $\overline{x}$ and $-\overline{x}$. The topology on $\mathbb{P}_2$ is Euclidean one, with base $\mathcal{D}_{\mathrm{pr}}$. We transfer
this base via $F$ to $\mathbb{P}_2'$, but we will only use the Euclidean topology on $R$. We define the {\em projective lines in $\mathbb{P}_2'$} as 
the intersections
$$
L_{\mathrm{pr}}=H\cap\mathbb{P}_2'
$$
where $H$ is a~plane in $\R^3$ going through the origin. The finite part of $L_{\mathrm{pr}}$ in $R$ is an ordinary affine line or is empty.
Clearly, $F$ maps every projective line to an $L_{\mathrm{pr}}$, and $F^{-1}$ maps every $L_{\mathrm{pr}}$ to a~projective line. In the proof 
of Theorem~\ref{proj_hilb_thm} we use yet another representation of the real projective plane, namely as the set of equivalence classes 
$$
(\R^3\setminus\{(0,0,0)\})/\!\sim\;,
$$ 
where two triples are equivalent in $\sim$ if one is a~nonzero multiple of the other. We write 
$(x:y:z)$ for the representatives of the equivalence classes, to emphasize that only mutual ratios of the entries in the triple are relevant for determination of the equivalence
class, the projective point.

We define the sets $\mathcal{P}_{\mathrm{pr}}$, $\mathcal{L}_{\mathrm{pr}}$, $\mathcal{C}_{\mathrm{pr}}$, and $\mathcal{D}_{\mathrm{pr}}$ as consisting 
of all projective points, all projective lines, all projective circles, and all projective discs, respectively. By $\mathcal{S}_{\mathrm{pr}}$ we 
denote a~possibly empty set system of nonempty subsets of the projective plane. The first four sets are pairwise disjoint and we assume that also
$\mathcal{S}_{\mathrm{pr}}$ is disjoint to each of the first three sets. A~{\em projective Euclidean construction with  $\mathcal{S}_{\mathrm{pr}}$-arbitrary
points}, abbreviated $\mathrm{EC}_{\mathrm{pr}}(\mathcal{S_{\mathrm{pr}}})$, is a~rooted tree $T=(r,V,E)$ with the vertices
$$
V\sus(\mathcal{P}_{\mathrm{pr}}\cup\mathcal{L}_{\mathrm{pr}}\cup\mathcal{C}_{\mathrm{pr}}\cup\mathcal{S}_{\mathrm{pr}})^*,\;\text{ with the root }\;
r\in(\mathcal{P}_{\mathrm{pr}}\cup\mathcal{L}_{\mathrm{pr}}\cup\mathcal{C}_{\mathrm{pr}})^*\cap V\;,
$$
and with the parent$-$child pairs determined by modifications of the six rules in Definition~\ref{defi_EC}. In more details, if $u=a_1a_2\ds a_m\in V$,
$m\in\N_0$, is a vertex of $T$, its child $v=a_1a_2\ds a_{m+1}$ is determined by exactly one of the following six rules, where in the first five we assume 
that $a_m\not\in\mathcal{S}_{\mathrm{pr}}$. 1. There is no $v$ and $u$ is a~leaf of $T$. 2. The vertex $v$ is the only child of $u$ and $a_{m+1}$ is 
a~projective line determined by two distinct projective points in $u$. 3. The vertex $v$ is the only child of $u$ and $a_{m+1}$ is a~projective circle 
or a~projective point $k(p,q,r)$ determined by three (not necessarily distinct) projective points $p$, $q$ and $r$ in $u$. 4. The vertex $v$ is the only 
child of $u$ and $a_{m+1}$ is an intersection projective point of two distinct projective lines, or two distinct projective circles, or a~projective line
and a~projective circle in $u$. 5. The vertex $v$ is the only child of $u$ and $a_{m+1}\in\mathcal{S}_{\mathrm{pr}}$. 6. We have 
$a_m\in\mathcal{S}_{\mathrm{pr}}$ and for every projective point $a_{m+1}\in a_m$, the vertex $v=a_1a_2\ds a_{m+1}$ is a~child of $u$. 

Further notions pertaining to $\mathrm{EC}(\mathcal{S})$ $T=(r,V,E)$, namely branches in $T$, the type of $T$ (compass, straightedge and general), 
terminating $T$, classical $T$, deterministic $T$, and $T$ constructing $\mathcal{K}$ from $r$, are adapted to 
$\mathrm{EC}_{\mathrm{pr}}(\mathcal{S_{\mathrm{pr}}})$ straightforwardly and we skip details. It would be also straightforward to define the variant 
 $\mathrm{EC}_{\mathrm{pr},\infty}(\mathcal{S_{\mathrm{pr}}})$ of  $\mathrm{EC}_{\mathrm{pr}}(\mathcal{S_{\mathrm{pr}}})$ with infinite starting 
 configuration.

We proceed to projective version of Hilbert's theorem.

\begin{thm}[projective Hilbert's theorem]\tec\label{proj_hilb_thm}
Let $k\sus\mathbb{P}_2$ be a~projective circle, $c\in\mathbb{P}_2$ be its center, and let $\mathcal{K}=\{c\}$. 
Then every straightedge classical (not necessarily terminating) $\mathrm{EC}_{\mathrm{pr}}(\mathcal{D_{\mathrm{pr}}})$ $T=(r,V,E)$ 
with $r=k$ has a~branch $B$ such that $c\not\in u$ for every vertex $u\in B$. So no  straightedge classical $\mathrm{EC}_{\mathrm{pr}}(\mathcal{D_{\mathrm{pr}}})$ constructs $\mathcal{K}$ (the center $c$ of $k$) from $r$ (the given projective circle $k$).
\end{thm}
\proof
Let $k$ and $c$ be as stated. Like in affine Hilbert's theorem, we need a~set $X\sus\mathbb{P}_2$ of projective points with the next properties.
\begin{enumerate}
    \item $X$ is dense in $\mathbb{P}_2$, which means that $X\cap D\ne\emptyset$ for any projective disc $D$. 
    \item $c\not\in X$.
    \item $X$ is {\em $\mathrm{H}_{\mathrm{pr}}$-closed}, any two projective lines determined by two pairs of distinct projective points in $X$
    have intersection in $X$, and any projective line determined by two distinct projective points in $X$ either misses $k$ or intersects $k$ in one 
    or two projective points in $X$.
\end{enumerate}
Assuming that such set $X$ exists, we proceed as before and inductively define in any given $\mathrm{EC}_{\mathrm{pr}}(\mathcal{D}_{\mathrm{pr}})$ a~branch 
$B$ such that in every vertex $u\in B$ every projective point letter lies in $X$ and thus differs from $c$. It remains to show that $X$ exists.

As before we find $X$ first for a particular projective circle $k$ and then get it by transformations for any given projective circle $k'$.  We let 
the $k$ be
$$
k=\{\{\overline{x},\,-\overline{x}\}\in\mathbb{P}_2\;|\;\text{$\overline{x}=(x,\,y,\,z)$ satisfies $x^2+y^2=1/2$}\}\;.
$$
Thus $k$ is formed by the intersections of $\mathbb{S}$ with the lines going through the origin and making with the plane $z=0$ angle $\pi/4$. The 
center of $k$ is $c=\{(0,0,1),(0,0,-1)\}$, the poles of $\mathbb{S}$. The above described bijection $F\cc\mathbb{P}_2\to\mathbb{P}_2'$ sends
the projective circle $k$ to the circle $k_0\sus(z=1)$ in the plane $R$. The circle $k_0$ has center $c_0=(0,0,1)$ and radius $1$. 

We lift Strommer's partial map in the plane $R=(z=1)$ as $f(x,y,1)$ where $f(x,y,z)=(\frac{1}{x},\frac{y}{x},z)$. It fixes the circles 
$$
k(a):=\{(x,y,1)\;|\;(x-a)^2+y^2=a^2-1\}\sus R,\ a>1\;. 
$$
We get an analogous map $f_{\mathrm{pr}}\cc\mathbb{P}_2\to\mathbb{P}_2$ fixing $k$. First we set $a=\sqrt{2}$, conjugate $f$ by the shift 
$\sigma(x,y,z)=(x-\sqrt{2},y,z)$, and get a~partial map $f_0$ fixing the circle $k_0$. Then we extend $f_0$ to everywhere defined map $\overline{f_0}\cc\mathbb{P}_2'\to\mathbb{P}_2'$. Finally we conjugate $\overline{f_0}$ by $F$. 

Thus
$$
f_0(x,\,y,\,z)=\sigma\circ f\circ\sigma^{-1}
=\bigg(\frac{1}{x+\sqrt{2}}-\sqrt{2},\,\frac{y}{x+\sqrt{2}},\,z\bigg)
$$
and
$$
\text{$f_0(x,\,y,\,1)\cc R\setminus L_0\to R\setminus L_0$, for the line $L_0=(-\sqrt{2},\,y,\,1)\sus R$}\;.
$$
Clearly, $f_0(x,y,1)$ is continuous. We extend it to $L_0\cup L\cup P$ by the projectivization 
$$
\overline{f_0}(x:y:z):=f_0(x/z,\,y/z,\,1)=\left(-\sqrt{2}x-z:y:x+\sqrt{2}z\right)\cc\mathbb{P}_2'\to\mathbb{P}_2'\;.
$$
Clearly, the value of $\overline{f_0}$ depends only on $(x:y:z)$, and $\overline{f_0}$ coincides with $f_0(x,y,1)$ 
on $R\setminus L_0$. Out of thin air we got the new values 
$$
\text{$\overline{f_0}(-\sqrt{2}:y:1)=(1/y:1:0)\in L$ for $y\ne0$, $\overline{f_0}(-\sqrt{2}:0:1)=(1:0:0)\in P$}\;,
$$
$\overline{f_0}(0:1:0)=(0:1:0)\in L$ (a~fixed point of $\overline{f_0}$), and $\overline{f_0}=\overline{f_0}^{-1}$ on $L\cup P$. So $\overline{f_0}$ is an involution.
The values of $\overline{f_0}$ are linearly independent homogeneous linear polynomials, and therefore $\overline{f_0}$
maps every projective line in the projective plane $(x:y:z)$ to another such line. Hence $\overline{f_0}$ maps every projective line in $\mathbb{P}_2'$
to another such line. Note that also
$$
\overline{f_0}(k_0)=f_0(k_0)=(\sigma\circ f\circ\sigma^{-1})(k_0)=(\sigma\circ f)(k(\sqrt{2}))=\sigma(k(\sqrt{2}))=k_0\;.
$$
We finally define
$$
f_{\mathrm{pr}}:=F^{-1}\circ\overline{f_0}\circ F\cc\mathbb{P}_2\to\mathbb{P}_2\;.
$$

Thus $f_{\mathrm{pr}}$ is an involution and hence a~bijection. It is also clear that $f_{\mathrm{pr}}$ is continuous
on
$$
\mathbb{P}_2\setminus(F^{-1}(L_0)\cup E)\;\text{ where }\;E=\{\{(x,\,y,\,0),\,(-x,\,-y,\,0)\}\in\mathbb{P}_2\;|\;x,\,y\in\R\}\;,
$$
on $\mathbb{P}_2$ with two projective lines deleted.
The map $f_{\mathrm{pr}}$ fixes the projective circle $k$ because 
$$
f_{\mathrm{pr}}(k)=(F^{-1}\circ\overline{f_0}\circ F)(k)=(F^{-1}\circ\overline{f_0})(k_0)=F^{-1}(k_0)=k\;.
$$
As for the center $c=\{(0,0,1),(0,0,-1)\}$,
\begin{eqnarray*}
f_{\mathrm{pr}}(c)&=&(F^{-1}\circ\overline{f_0}\circ F)(c)=(F^{-1}\circ\overline{f_0})(0,\,0,\,1)=F^{-1}(-1/\sqrt{2},\,0,\,1)\\
&=&\{(-1/\sqrt{3},\,0,\,\sqrt{2/3}),\,(1/\sqrt{2},\,0,\,-\sqrt{2/3})\}=:c'\ne c\;.
\end{eqnarray*}
Finally, $f_{\mathrm{pr}}$ maps any projective line to a~projective line. This is immediate from the fact that $\overline{f_0}$ 
preserves the projective lines in $\mathbb{P}_2'$ and from their definition.

A~single map $f_{\mathrm{pr}}$ with these properties does not suffice, we need uncountably many of them. We get them as in the proof of Proposition~\ref{map_fp}.
Let $k_1$ be the projective circle centered at $c$
and going through the projective point $c'=f_{\mathrm{pr}}(c)$, let $p\in k_1$, and let $\tau_p$ be the rotation of $\mathbb{P}_2$ around the $z$-axis 
(i.e., the line $(0,0,z)$) that moves $c'$ to $p$. We set
$$
f_{\mathrm{pr},p}=\tau_p\circ f_{\mathrm{pr}}\;.
$$
It is easy to see that every $f_{\mathrm{pr},p}$ shares with $f_{\mathrm{pr}}=f_{\mathrm{pr},c'}$ the above properties: it is a~bijection, 
is continuous except possibly on two projective lines, preserves the projective circle $k$, and preserves the set of projective lines. Also, $f_{\mathrm{pr},p}(c)=p$. 

We take any countable set $X'\sus\mathbb{P}_2$ that is dense in $\mathbb{P}_2$ and is $\mathrm{H_{pr}}$-closed. Then we set 
$$
\text{$X=f_{\mathrm{pr},p}^{-1}(X')$, for any projective point $p\in k_1\setminus X'$}\;.
$$
It follows that $X$ has the properties 1--3 and is the desired set.

As for the general case, suppose that $k'$ is any given projective circle and $c'$ is its center. We argue somewhat differently compared to the 
affine case because we do not see what could be an analog of the similarity map with a~given center for the projective plane. But we can begin this proof
with the projective circle $k$ with the same center $c=\{(0,0,1),(0,0,-1)\}$ but with the radius equal to that of $k'$, we simply take at the start 
an appropriate parameter $a>1$ and shift $\sigma(x,y,z)=(x-a,y,z)$. Then we rotate $\mathbb{P}_2$ around an appropriate axis going through 
the origin so that $c$ moves to $c'$
and hence $k$ to $k'$. If we denote this rotation by $g$, and by $X\sus\mathbb{P}_2$ the set with properties 1--3 with respect to the modified projective 
circle $k$, then the set $g(X)$ has properties 1--3 with respect to the projective circle $k'$.
\eproof

\noindent
We conclude with some remarks. To our knowledge, the previous theorem is the first treatment of Hilbert's theorem in projective form.
One can shorten the proof at the cost of omitting the affine motivation and start directly from the involution
$$
\overline{f_0}(x:y:z)=\left(-\sqrt{2}x-z:y:x+\sqrt{2}z\right)
$$
of the projective plane $(x:y:z)$. Then the model $\mathbb{P}_2'$ is not needed and it is easy to check that $\overline{f_0}$ preserves the 
projective circle $k$. Theorem~\ref{proj_hilb_thm} in its present form does not directly imply Theorem~\ref{hilb_thm} because the 
map $F$ does not in general send projective circles to circles in $R$.

\end{document}